\documentclass[10pt]{article}
\usepackage{amsfonts}
\usepackage{amsmath,amssymb}
\usepackage[a4paper]{geometry}

\newtheorem{Lemma}{Lemma}[section]
\newtheorem{theorem}[Lemma]{Theorem}
\newtheorem{proposition}[Lemma]{Proposition}
\newtheorem{lemma}[Lemma]{Lemma}

\newtheorem{definition}[Lemma]{Definition}

\usepackage{xcolor}

\begin{document}

\begin{center}
{\Large \textbf{Well-posedness of stochastic third grade fluid equation}} \vspace{%
3mm}\\[0pt]

\vspace{1 cm}

\textsc{Fernanda Cipriano}{\footnote{%
Departamento de Matem\'atica, Faculdade de Ci\^encias e Tecnologia da
Universidade NOVA de Lisboa and Centro de Matem\'atica e Aplica\c c\~oes.
E-mail: cipriano{\char'100}fct.unl.pt.},}
\hspace{0mm}
\textsc{\ Philippe Didier}{\footnote{%
Departamento de Matem\'atica, Faculdade de Ci\^encias e Tecnologia da
Universidade NOVA de Lisboa and Centro de Matem\'atica e Aplica\c c\~oes. E-mail:  pld{\char'100}fct.unl.pt.},} 
\hspace{0mm}
\textsc{\ S\'{i}lvia Guerra}{\footnote{%
Universidade NOVA de Lisboa and Centro de Matem\'atica e Aplica\c c\~oes, 
 and NOVA School of Business and Economics. E-mail:  si.guerra{\char'100}campus.fct.unl.pt.},} 

\end{center}

\noindent 

\date{\today }

\begin{abstract}
In this paper, we establish the well-posedness for the third grade fluid equation perturbed by a multiplicative white noise.
This equation describes the motion  of a non-Newtonian fluid of differential type with relevant viscoelastic properties.  We are faced with a strongly nonlinear stochastic partial differential equation supplemented with a Navier slip boundary condition.
Taking the initial condition in the Sobolev space $H^2$, we show the existence and the uniqueness  of the strong (in the probability sense) solution in a two dimensional and  non axisymmetric bounded domain. 

\vspace{3mm}

\vspace{3mm}\newline


\vspace{3mm}
\end{abstract}

\section{Introduction}

\setcounter{equation}{0}
This work deals  with fluids of differential type, which is a very special class of non-Newtonian fluids. For these fluids,  
the relation between the shear stress  and the shear strain rate is not linear, this means that  the viscosity does not satisfy   Newton's law.
Considering  the velocity field   $y$ of the fluid, we introduce  the Rivlin-Ericksen kinematic  tensors in \cite{RE55} ${\bf A}_n,$ $n\geq 1$,   defined by 
\begin{align}
{\bf A}_1( y)&=
	\nabla  y+\nabla  y^\top,\\
		{\bf A}_n(y)&=\frac{d}{dt}{\bf A}_{n-1}(y)+{\bf A}_{n-1}(\nabla y)+(\nabla y)^\top {\bf A}_{n-1}(y), \quad n=1,2,\dots,
	\end{align}
While a Newtonian fluid  is characterized by  the relation $\mathcal{T}=\nu {\bf A}_1( y)$, 	where $\mathcal{T}$ is the shear  stress tensor and $\nu$ is the viscosity;	
a third grade fluid is given  by the following constitutive equation 
\begin{equation}\label{cauchy}
	\mathbb{T}=-p\mathbb{I}+\nu {\bf A}_1(y)+\mathbb{T}_1+\mathbb{T}_2,
	\end{equation}
where $\mathbb{T}$ is the Cauchy stress tensor,
 $\frac{D}{Dt}=\frac{d}{dt}+ y\cdot \nabla$ stands for the 
material derivative, $p$ is the pressure,  $\mathbb{I}$ is the unit tensor, 
\begin{align*}
\mathbb{T}_1&=\alpha_1 {\bf A}_2+\alpha_2{\bf A}_1^2,\\
\mathbb{T}_2&=\beta_1{\bf A}_3+\beta_2({\bf A}_1{\bf A}_2+{\bf A}_2{\bf A}_1)+\beta_3({\rm tr}{\bf A}_1^2){\bf A}_1,
\end{align*}
and  $\alpha_1$,$\alpha_2$, $\beta_1$, $\beta_2$ and $\beta_3$ are material moduli.
The momentum equations	is given by
$$
\frac{Dy}{Dt}={\rm div} \,\mathbb{T}.
$$
According to the analysis in  \cite{DR95}, \cite{FR80}
in order to allow the motion of the fluid to be  compatible with thermodynamic, it  should be imposed that
\begin{equation}
\label{10:30_03}
	\nu\geq 0, \qquad
	\alpha_1\geq 0,\qquad |\alpha_1+
	\alpha_2|\leq \sqrt{24\nu\beta_3}, \qquad \beta_1=\beta_2=0, \qquad \beta_3\geq 0.
	\end{equation}
Setting 
\begin{equation}
\label{10:31_03}
	\beta=\beta_3, \qquad  A={\bf A}_1,
\end{equation}
	and denoting $\upsilon (y)=y-\alpha_1\Delta y$, the deterministic  equation of motion for a third grade fluid reads
	\begin{equation}
\begin{array}{cccc}
 \displaystyle \frac{d}{dt}(\upsilon (y))=-\nabla p+
 \nu \Delta y-
(Y\cdot \nabla)\upsilon-\sum_j\upsilon^j\nabla y^j+(\alpha_1+\alpha_2)\,{\rm div}
\left(A^2\right) 
 +\beta{\rm div}\left(|A|^2A\right) +U,
 \end{array}
\label{equation_etat}
\end{equation}
where $U$ represents a body force.
In practice, the non-Newtonian fluids are real fluids present in  industry, food processing, biological fluids, etc. (see 
\cite{FR80}, \cite{DR95}, 
\cite{HKM03}, \cite{HSA07}, \cite{RKWA17}). In the case where $\beta=0$, the model reduces to the 
 second grade fluid model which is  mathematically more tractable. However the second grade fluid model  does not capture  important rheological
properties as for instance the  shear thinning and shear thickening effects, so there is  
a real need to study the  third grade fluid model.

From a theoretical   point of view,   this model is  a strongly nonlinear partial differential equation modelling  complex viscoelastic fluids, therefore the small noise perturbations 
should have relevant impact in the fluid dynamic. It is well known that increasing  the typical velocity will increase the  Reynolds number:  the fluid develops a turbulent behaviour and  small disturbances should have strong macroscopic effects on the dynamic. 
Let us mention the pioneer work \cite{B95} on the stochastic Navier-Stokes equations, and  \cite{BT73} (see Lemma 2.2) where 
the stochastic Navier-Stokes equations were deduced from fundamental principles 
 showing that the stochastic  equations are  real physical
models. There is an extensive literature on the stochastic Eulerian description of a Newtonian fluid
(see also \cite{CC07}, \cite{CP19} and references therein for a stochastic Lagrangian approach). 
In this paper, we perturb the deterministic third grade fluid model by a multiplicative white noise 
\begin{align}
\label{11}
d(\upsilon (Y))&=\bigl(-\nabla p+
 \nu \Delta Y-
(Y\cdot \nabla)\upsilon-\sum_j\upsilon^j\nabla Y^j+(\alpha_1+\alpha_2)\,{\rm div}
\left(A^2\right)\notag\vspace{2mm} \\
\vspace{2mm} 
 &\qquad+\beta{\rm div}\left(|A|^2A\bigr) +U\right)\,dt
 + \sigma(t,Y)\,d\mathcal{W}_{t},
\end{align}
and we extend the deterministic results on the  existence and  uniqueness obtained in
 \cite{BI04} and  \cite{BI06}, where the parameter $\beta$ is considered strictly positive,  the equation is supplemented with a Navier slip boundary condition and the  initial condition is taken in  the Sobolev space $H^2$. As far as we know, the  stochastic third grade fluid equations are  being studied  for the first time in this work. Here, we also consider 
$\beta>0$ and the initial condition in $H^2$.

Concerning the  analysis, we recall  that the deterministic strategy in  \cite{BI04} and  \cite{BI06} is based on the deduction of a priori estimates that allow to use compactness theorems in order to pass the non linear terms to the limit in the weak sense.
 However, for the stochastic equation, due to the lack of regularity with respect to time and to the stochastic parameter,  we can not use the compactness arguments to pass to the limit in those nonlinear terms. 
 Instead, we will follow the methods introduced in \cite{B99}, which have been  successfully applied to 
  the stochastic second grade fluids in \cite{CC16},    \cite{RS12} (see also \cite{CC18}, 
  \cite{C19}, \cite{CP19}, \cite{HAKAWA18}, \cite{R18}, \cite{SJT19}, \cite{ZZZ18}).  
  More precisely, we consider an appropriate Galerkin basis, and deduce suitable uniform estimates in order to get weak convergence of a subsequence. Next we project the weak limit on the finite $n-$dimensional space. we then show that the difference between 
the sequence obtained by projecting the weak limit and the finite dimensional Galerkin approximations converges strongly to zero up to a special stopping time. Finally we are able to identify the nonlinear terms of the equation.

The remainder of this paper is structured into  four sections and one appendix.
 Section 2 is devoted to the introduction of the 
functional settings and appropriate  notations.
 In Section 3 we present some  properties of the   nonlinear terms
of equation (\ref{11}) which will be applied in the next section.
 The main results concerning the existence and uniqueness of the strong stochastic
solution is established in Section 4. 
Finally, in the Appendix we collect some known inequalities related with the nonlinear terms of the equation that are used throughout the article.

\section{Functional setting and notations}

\setcounter{equation}{0}
We consider the  stochastic third  grade fluid equation \eqref{11} in a bounded, not axisymmetric and simply
connected domain $\mathcal{O}$ of $\mathbb{R}^{2}$ with a sufficiently
regular boundary $\Gamma $, and supplemented with a Navier slip boundary condition, which reads 
\begin{equation}
\left\{
\begin{array}{cccc}
d(\upsilon (Y))=
 \bigl(-\nabla p+\nu \Delta Y-
(Y\cdot \nabla)\upsilon-\sum_j\upsilon^j\nabla Y^j+(\alpha_1+\alpha_2)\,{\rm div}
\left(A^2\right) \vspace{2mm} \\
\vspace{2mm} 
 +\beta{\rm div}\left(|A|^2A\right) +U \bigr)dt
 + \sigma(t,Y)  \,d\mathcal{W}_{t},
 &  \multicolumn{1}{l}{\mbox{in}\
\mathcal{O}\times (0,T),} \vspace{2mm}\\
\multicolumn{1}{l}{\mathrm{div}\,Y=0} & \multicolumn{1}{l}{\mbox{in}\
\mathcal{O}\times (0,T),} \vspace{2mm}\\
\multicolumn{1}{l}{Y\cdot \mathrm{n}=0,\qquad
 (\mathrm{n}\cdot
D(Y))\cdot \mathrm{\tau }=0} & \multicolumn{1}{l}{\mbox{on}\
\Gamma \times (0,T),}\vspace{2mm} \\
\multicolumn{1}{l}{Y(0)=Y_{0}} & \multicolumn{1}{l}{\mbox{in}\ \mathcal{O},}
\end{array}
\right.  \label{equation_etat}
\end{equation}
\vspace{1mm}\\
where $Y$ is the velocity field of the fluid, $\nabla Y$ is its  Jacobian matrix, $D(Y)=\frac{\nabla Y+(\nabla Y)^\top}{2}$,
$A=A(Y)=2D(Y)$,  
$\upsilon (Y)=Y-\alpha_1\Delta Y$ and the
 constants $\nu$, $\alpha_1$, $\alpha_2$, $\beta$ verify \eqref{10:30_03}-\eqref{10:31_03}.
 The  stochastic
perturbation is defined by
$$\sigma(t,Y) \,d\mathcal{W}_{t}=\sum_{k=1}^m  \sigma^k(t,Y) \,d\mathcal{W}^k_{t},$$ 
 where the diffusion coefficient 
$$\sigma(t,Y)=(\sigma^1(t,Y), \dots, \sigma^m (t,Y))$$
satisfy
suitable growth assumptions that will be defined below, and  $\mathcal{W}_{t}
=(\mathcal{W}^1_{t}, \dots, \mathcal{W}^m_{t})$ is a standard $\mathbb{R}^{m}$-valued Wiener
process defined on a complete probability space $(\Omega ,\mathcal{F},P)$
endowed with a filtration $\left\{ \mathcal{F}_{t}\right\} _{t\in \lbrack
0,T]}$ for $\mathcal{W}_{t}$. We assume that $\mathcal{F}_{0}$\ contains every $P$-null subset of $%
\Omega $.

We consider the following Hilbert spaces
\begin{equation}
\begin{array}{l}
H=\left\{ y\in L^{2}(\mathcal{O})\mid \mathrm{div}\,y=0\ \text{ in }\mathcal{%
O}\ \mbox{ and }\ y\cdot \mathrm{n}=0\ \mbox{ on }\Gamma \right\} ,\vspace{%
2mm} \\
V=\left\{ y\in H^{1}(\mathcal{O})\mid \mathrm{div}\,y=0\ \mbox{ in }\
\mathcal{O}\mbox{ and }\ y\cdot \mathrm{n}=0\ \text{ on }\ \Gamma \right\} ,%
\vspace{2mm} \\
W=\left\{ y\in V\cap H^{2}(\mathcal{O})\mid (\mathrm{n}\cdot
D(y)) \cdot \mathrm{\tau }=0\ \ \mbox{on}\ \Gamma \right\}.
\end{array}
\label{w}
\end{equation}
On $H$ we consider the $L^2-$inner product $(\cdot,\cdot)$ and the associated norm $\|\cdot\|_2.$ 

Let us introduce the Helmholtz projector $\mathbb{P}:L^{2}(\mathcal{O}%
)\longrightarrow H$, which  is the  linear bounded
operator  characterized  by the following $L^2-$orthogonal decomposition
\begin{equation*}
v=\mathbb{P}v+\nabla \phi ,\qquad  \phi
\in H^{1}(\mathcal{O}).
\end{equation*}

We define the following inner products
\begin{eqnarray}
\label{ip_V}
\left( u,z\right) _{V} &:=&\left( \upsilon (u),z\right) =\left( u,z\right)
+2\alpha_1 \left( Du,Dz\right), \\
\left( u,z\right) _{W} &:=&\left( u,z\right) _{V}+\left( \mathbb{P}\upsilon (u),\mathbb{P}\upsilon
(z)\right) ,
\label{ip_W} 
\end{eqnarray}
and denote by $\|\cdot\|_V$ and $\|\cdot\|_W$ the norms induced by these inner products on $V$ and $W$, respectively.

The space $V$ being a subspace of $H^1(\mathcal{O})$ is naturally endowed with the  norm $\|\cdot\|_{H^1}$.
The  norms $\|\cdot\|_{H^1}$ and 
$\|\cdot\|_V$ on the space $V$ are equivalent. 
Similarly  $W\subset H^2(\mathcal{O})$ and the  norms  $\|\cdot\|_W$ and $\|\cdot\|_{H^2}$ are equivalent on $W$.

In addition, through the article the usual norms on
the spaces  $L^p(\mathcal{O})$, $p>0$, 
are denoted by
 by $\|\cdot\|_p$ and the norms on the Sobolev spaces $W^{1,p}(\mathcal{O})$ for $p>2$  are denoted  by 
$\|\cdot\|_{W^{1,p}}.$

\bigskip

Let us introduce   the  trilinear functional 
\begin{equation}
\label{tf_1}
b(\phi ,z,y)=\left( \phi \cdot \nabla z,y\right) ,\qquad \forall \phi
,z,y\in V.
\end{equation}
Taking into account that
 $\phi $ is divergence
free and $(\phi \cdot n)=0$ on $\Gamma $, a standard integration by parts gives 
\begin{equation}
b(\phi ,z,y)=-b(\phi ,y,z).  \label{ra}
\end{equation}%

\bigskip

Assume that $\sigma=\left(\sigma^1, \dots,\sigma^m\right): [0,T]\times V\rightarrow  (L^2(\mathcal{O}))^{m} $ is   Lipschitz in the second variable and verifies a growth condition, i.e.,
 there exist 
 positive constants $L,\,K$ and $0\leq \gamma <2$ such that
\begin{align}
\label{LG_1}
\left\Vert \sigma(t,y)\right\Vert _2^2& \leq L(1+\|y\|^\gamma_{W^{1,4}})
,\quad \forall y\in W^{1,4}(\mathcal{O})\cap V, \\
\left\Vert \sigma(t,y)-\sigma(t,z)\right\Vert _{2}^{2}& \leq K\left\Vert
y-z\right\Vert _{V}^{2} ,\qquad \forall y,z\in V,\;t\in \lbrack 0,T],
 \label{LG_2}
\end{align}
where
\begin{equation*}
\left\Vert \sigma(t,y)\right\Vert _{2}^2:=\sum_{i=1}^{m}\left\Vert
\sigma^{i}(t,y)\right\Vert _{2}^2 .
\end{equation*}
We also introduce the notation
\begin{equation*}
|\left( \sigma(t,y),v\right) | : = \left( \sum_{k=1}^{m}\left(  \sigma^k(t,y), v \right)^2 \right)^{1/2}, \quad \forall v\in L^2(\mathcal{O}).  
\end{equation*} 
In addition, we take $p\geq 6$
and suppose that the initial condition $Y_0$ and the force $U$ satisfy \begin{equation}
\label{R_1}
Y_0\in L^p(\Omega,W), \quad \text{there exists }\lambda>0\text{ such that}\quad
\mathbb{E}\text{e}^{\lambda\bigl(\int_0^T\|U\|_2^2ds+\|Y_0\|^2_V\bigr)} <\infty.
\end{equation} 

\bigskip

Here, we   recall the  
 Korn inequality
\begin{equation}
\label{K_1}
\|y\|_{W^{1,p}}\leq K_1(p)\left(\|y\|_p+\|A(y)\|_{p}\right), \quad \forall y\in V, \quad p\geq 2,
\end{equation}
and the Poincar\'e inequality
\begin{equation}
\label{Poincare}
\|y\|_{2}\leq \mathcal{P}\|\nabla y\|_2, \quad \forall y\in V.
\end{equation}

For  non axisymmetric bounded domains, we also have the following version of the  Korn inequality (see Theorem 3 in \cite{DV02})

\begin{equation}
\label{14:37_3Ag}
\|\nabla y\|_2\leq K_2(\mathcal{O})\|A(y)\|_{2}, \quad \forall y\in V.
\end{equation}
The Sobolev 
 embedding $H^1(\mathcal{O})\hookrightarrow L^4(\mathcal{O})  $
 and \eqref{Poincare} give
\begin{equation*}
\|y\|_{4}\leq K_3\|\nabla y\|_{2}, \quad \forall y\in V.
\end{equation*}
Combining this inequality  with \eqref{14:37_3Ag}, we get
\begin{equation}
\label{14:45_04:08}
\|y\|_{4}\leq K_3\|\nabla y\|_{2}\leq K_3 K_2(\mathcal{O})\|A(y)\|_{2}. 
\end{equation}
Due to the  embedding $L^4(\mathcal{O}) \hookrightarrow L^2(\mathcal{O})$, we have 
\begin{equation}
\label{14:27_04:08}
\|y\|_2\leq C_*\|y\|_4.
\end{equation}
 Then \eqref{14:45_04:08}, \eqref{14:27_04:08} and \eqref{K_1} yield the following   lemma:
\begin{lemma}
\label{K_2} There exists  a positive constant  $K_* $ 
such that
\begin{equation}
\label{Korn_*}
\|y\|_{W^{1,4}}\leq K_* \|A(y)\|_{4}, \quad \forall y\in V.
\end{equation}
\end{lemma}

\bigskip

Let us mention that through the article, we represent by $C$ a generic constant. Its value can change  from line to line.
To explicitly write its dependence with respect  of some  parameters $\lambda_1,\dots,\lambda_k$, we also write 
$C(\lambda_1,\dots,\lambda_k)$ instead of $C$. 

\bigskip

We end this section with the Young's inequality 
\begin{equation}
\label{R_15}
uz\leq\frac{1}{r}u^r+\frac{1}{s}z^s,\quad  \forall u,z\geq 0, \quad s,r>0\quad\text{such that}\quad \frac{1}{r}+\frac{1}{s}=1.
\end{equation}
Accordingly, for real numbers $\gamma, a, b, x$ such that $0\leq \gamma<a$ and $b,x\geq 0$, we  have the  algebraic relation
\begin{equation}
\label{R_16}
\forall \delta>0,\qquad bx^\gamma\leq C(\gamma, a,b,\delta)+\delta x^a,
\end{equation}
that will be used several times.

\section{Preliminary   results}
\setcounter{equation}{0}

We consider the following auxiliary modified Stokes problem with Navier boundary
condition
\begin{equation}
\left\{
\begin{array}{ll}
\tilde{f}-\alpha_1 \Delta \tilde{f}=f-\nabla p,\qquad \mathrm{div}\,\tilde{f}=0 & \quad \mbox{in}\
\mathcal{O},\vspace{2mm} \\
\tilde{f}\cdot \mathrm{n}=0,\qquad  (\mathrm{n}\cdot D(\tilde{f}))
\cdot \mathrm{\tau }=0 & \quad \mbox{on}\ \Gamma \vspace{2mm}.%
\end{array}%
\right.  \label{GS_NS}
\end{equation}
We  recall from  \cite{CC18} that assuming  $f\in H^{m}(\mathcal{O})$, $m=0,1$, the problem
\eqref{GS_NS} has a solution $(\tilde{f},p)\in H^{m+2}(\mathcal{O})\times H^{m+1}(%
\mathcal{O})$ verifying 
\begin{align}
\Vert \tilde{f}\Vert _{H^{2}}& \leq C\Vert f\Vert _{2}.  \label{SP_L2} 
\end{align}
According to the definition of the inner product \eqref{ip_V}, we have 
\begin{equation}
\label{12-08-36}
(\tilde{f}, h)_V
=(f, h),\quad \forall h\in V.
\end{equation} 

In the next two lemmas, we establish   properties  
  of the nonlinear terms that will be useful in Section \ref{section4} to identify the  weak limits of the nonlinear terms of the equation.
Let us introduce the operators
\begin{eqnarray}
&&S(y):=\beta\left(|A(y)|^2A(y)\right), \quad \label{operator1}\\
&&N(y):=\alpha_1\left(y\cdot \nabla A(y)+(\nabla y)^\top\, A(y)+A(y)\, \nabla y\right)
-\alpha_2(A(y))^2.\label{operator2}
\end{eqnarray}

\begin{lemma} For any $y, \,\hat{y}, \,\phi\in W$, we have 
\begin{align}
\label{29_5_S1-S22}
\left|\langle {\rm div} (S(y)-S(\hat{y}), \phi\rangle\right|\leq 
C\|y\|^2_{W}\|y-\hat{y}\|_{V}\|\phi\|_{W}+
C\|\hat{y}\|_{W}\left\||A|^2-|\hat{A}|^2\right\|_2\|\phi\|_{W},
\end{align}
where $A=A(y)$ and $\hat{A}=A(\hat{y}).$
\end{lemma}
\textbf{Proof.}
Using the  H\"older inequality, and the Sobolev injections $H^1(\mathcal{O}))\hookrightarrow L^p(\mathcal{O})$ for $p<\infty$ and $H^2(\mathcal{O})\hookrightarrow L^\infty(\mathcal{O})$, we derive
\begin{align}
\label{13-17-13}
\left|\langle {\rm div} \left(S(y)-S(\hat{y})\right), \phi\rangle\right|
&= \left|\beta\int_{\mathcal{O}}(|A|^2A-|\hat{A}|^2\hat{A})\cdot\nabla\phi\right|\notag\\
&=\left|\beta\int_{\mathcal{O}}\left(|A|^2(A-\hat{A})
+\hat{A}(|A|^2-|\hat{A}|^2)\right)\cdot\nabla\phi\right|\notag\\
&\leq C\||A|^2\|_4\|A(y-\hat{y})\|_2\|\nabla\phi\|_4+
C\|\hat{A}\|_4\left\||A|^2-|\hat{A}|^2\right\|_2\|\nabla\phi\|_4
\notag\\
&\leq C\|y\|^2_{H^2}\|y-\hat{y}\|_{H^1}\|\phi\|_{H^2}+
C\|\hat{y}\|_{H^2}\left\||A|^2-|\hat{A}|^2\right\|_2\|\phi\|_{H^2}.
\end{align}

\bigskip\ $\hfill \hfill \blacksquare $

\begin{lemma}
\label{Est_N1-N2<} 
For any $y,\; \hat{y},\,\phi \in W$, the following inequality holds
\begin{align}
\label{13-15-22}
\left|\langle{\rm div}\left(N(\hat{y})- N(y)\right),\phi\rangle\right|
&\leq 
C\epsilon \left\|A(y-\hat{y})\sqrt{|A|^2+|\hat{A}|^2}\right\|_2\|\phi\|_V\notag\\
&
+C\|\hat{y}-y\|_V\left(\|y\|_W+\|\hat{y}\|_W\right)\|\phi\|_W.
\end{align}
where $A=A(y)$ and $\hat{A}=A(\hat{y})$.
\end{lemma}

\textbf{Proof.} Here we apply the same reasoning that is done in 
\cite{BI04} to show the property \eqref{25_2_N1-N2}.
\begin{align}
\label{13-15-23}
\langle{\rm div}\left(N(\hat{y})- N(y)\right),\phi\rangle
&=\langle N(\hat{y})- N(y),\nabla\phi\rangle=
\frac{1}{2}\langle N(\hat{y})- N(y),A(\phi)\rangle\notag\\
&=-\frac{\alpha_2}{2}\int_\mathcal{O}(A^2-\hat{A}^2)\cdot A(\phi)
-\frac{\alpha_1}{2}\int_\mathcal{O}(y\cdot\nabla A -\hat{y}\cdot \nabla \hat{A})\cdot A(\phi)\notag\\
&-\frac{\alpha_1}{2}\int_\mathcal{O}\left((\nabla y)^\top A
+A\nabla y-(\nabla\hat{y})^\top\hat{A}-\hat{A\nabla\hat{y}}\right)
\cdot A(\phi)=I_1+I_2+I_3.
\end{align}
\begin{align}
\label{13-15-24}
|I_1|\leq C \left\||A(y-\hat{y})|\sqrt{|A|^2+|\hat{A}|^2}
\right\|_2\|A(\phi)\|_2.
\end{align}
Next, we use the properties of the trilinear  form, as well as the H\"older inequality, and the Sobolev injections $H^1(\mathcal{O})\hookrightarrow L^4(\mathcal{O}) $ and $H^1(\mathcal{O})\hookrightarrow L^\infty(\mathcal{O})$ in order to deduce that
\begin{align}
\label{13-15-25}
|I_2|&\leq \left|b(y,A,A(\phi))+b(\hat{y},\hat{A},A(\phi))\right|
\notag\\
&\leq \left|b(y,A-\hat{A},A(\phi))\right|+\left|b(y-\hat{y},\hat{A},A(\phi))\right|\notag\\
&= \left|b(y,A(\phi),A-\hat{A})\right|+\left|b(y-\hat{y},\hat{A},A(\phi))\right|\notag\\
&\leq\|y\|_\infty\|A(\phi)\|_{H^1}\|A-\hat{A}\|_2+\|y-\hat{y}\|_4\|\hat{A}\|_{H^1}\|A(\phi)\|_{4}\notag\\
&\leq C\|y-\hat{y}\|_{H^1}\left(\|y\|_{H^2}+\|\hat{y}\|_{H^2}\right)\|\phi\|_{H^2}.
\end{align}
For $I_3$ we have the same estimate as for $I_1$, namely
\begin{align}
|I_3|\leq C \left\||A(y-\hat{y})|\sqrt{|A|^2+|\hat{A}|^2}
\right\|_2\|A(\phi)\|_2.
\label{13-15-26}
\end{align}

\bigskip\ $\hfill \hfill \blacksquare $

\section{Existence of strong solution}
\label{section4}

\setcounter{equation}{0}

\bigskip

This section establishes   the main  results of the article.
More precisely  the  solution of the equation is constructed through the finite dimensional Galerkin approximation method.
 We first deduce key uniform estimates for the finite dimensional approximations in order to get a weakly convergent sequence.
Next, with the help of a suitable stopping time, and using the structure of the equation,  we  improve the convergence results. Finally, with these  new convergence results, we will be able to identify the nonlinear terms of the equation.
\vspace{1mm}

Let us to introduce the notion of the solution.
\begin{definition}
\label{4.1} Let
$U\in L^{2}(\Omega \times (0,T),L^2(\mathcal{O}))$ and $Y_{0}\in L^{2}(\Omega ,W).$ Then 
a stochastic process $Y\in L^{2}(\Omega ,L^{\infty }(0,T;W))$ is
a strong solution of $(\ref{equation_etat})$, if the following equation holds
\begin{align}
\left( \upsilon (Y(t)),\phi \right)& = \int_{0}^{t}
\bigl[ -2\nu \left(
D(Y),D(\phi) \right) 
+\left((Y\cdot \nabla) \phi, \upsilon(Y)\right)-
\sum_j\left(\upsilon^j(Y)\nabla Y^j,\phi\right)
\bigr] \,ds  \notag
\\
&-\int_{0}^{t}\left((\alpha_1+\alpha_2)\,
\left(A^2\right)+\beta\left(|A|^2A\right),\nabla \phi\right)ds
\notag
\\
& +\left( \upsilon (Y(0)),\phi \right) +\int_{0}^{t}\left( U(s),\phi \right)
\,ds+\int_{0}^{t}\left(  \sigma(s,Y(s)),\phi \right) \,d\mathcal{W}_s 
\label{var_form_state}
\end{align}%
 $\quad \text{for a.e.} \;(\omega,t)\in\Omega\times[0,T]$ and for all $ \phi \in V$, where the stochastic integral is defined by
$$
\int_{0}^{t}\left( \sigma(s,Y(s)),\phi \right) \,d\mathcal{W}_s=\sum_{k=1}^m \int_{0}^{t} \left(  \sigma^k(s,Y(s)),\phi \right) \,d\mathcal{W}^k_s.
$$
\vspace{2mm}\newline
\end{definition}
Now we  state the main result.
\begin{theorem}
\label{the_1}
 Assume \eqref{LG_1}-\eqref{R_1}. Then there exists a unique solution
  $Y$ to equation $(\ref{equation_etat})$
which belongs to
\begin{equation*}
 L^p(\Omega, L^\infty(0,T;W)).
\end{equation*}
\end{theorem}

\bigskip

In order to  show the existence of the solution, we apply the  \textit{Galerkin's approximation method} for an appropriate  basis. We recall that the injection operator $I:W\hookrightarrow V$ being a compact operator guarantees the existence of a basis $\{e_{i}\}\subset W$ of eigenfunctions to the problem
\begin{equation}
\left( v,e_{i}\right) _{W}=\lambda _{i}\left( v,e_{i}\right) _{V},\qquad \forall
v\in W,\quad i\in \mathbb{N},  \label{mu}
\end{equation}%
which is  an orthonormal basis in $V$ and an  orthogonal basis in $W$. in addition the sequence
 $\{\lambda _{i}\}$  of the corresponding eigenvalues 
 fulfils the properties: $\lambda _{i}>0$, $\forall i\in \mathbb{N}$, and $\lambda _{i}\rightarrow \infty $ as $%
i\rightarrow \infty .$ Since the  ellipticity of the equation (\ref{mu}) increases the
regularity of their solutions (see \cite{BR03}), we  may consider $\{e_{i}\}\subset H^{4}.$  \vspace{2mm}

We consider the finite dimensional space  $W_{n}=\mathrm{span}\,\{e_{1},\ldots ,e_{n}\}$,
and  introduce the Faedo-Galerkin
approximation of the system (\ref{equation_etat}). Namely, we look for a solution 
to the following stochastic differential equation 
\begin{equation}
\label{y1}
\left\{
\begin{array}{l}
d\left( \upsilon_n  ,\phi \right) =\bigl( \nu \Delta Y_{n}-
(Y_{n}\cdot \nabla)\upsilon_n-\sum_j\upsilon_n^j\nabla Y_{n}^j+(\alpha_1+\alpha_2)\,{\rm div}
\left(A_n^2\right)\vspace{2mm} \\
\vspace{2mm} 
\qquad\qquad \qquad+\beta{\rm div}\left(|A_n|^2A_n\right)+U, \phi \bigr)dt
 +\left( \sigma(t,Y_{n}),\phi \right) \,d\mathcal{W}_{t},
\qquad  \forall \phi \in
W_{n},
\vspace{2mm} \\
Y_{n}(0)=Y_{n,0},
\end{array}
\right. 
\end{equation}
where  
\begin{equation*}
Y_{n}(t)=\sum_{j=1}^{n}c_{j}^{n}(t)e_{j}. 
\end{equation*}
Here $Y_{n,0}$ denotes the projection of the initial condition $Y_{0}$ onto
the space $W_{n}$,
$\upsilon_n=Y_{n}-\alpha_1\Delta Y_{n}$ and $A_n=\nabla Y_{n}+(\nabla Y_{n})^\top$.

Due to the relation \eqref{mu},  the sequence 
 $\{\widetilde{e}_{j}=\frac{1}{\sqrt{\lambda _{j}}}
e_{j}\}$ is an orthonormal basis for $W$ and
\begin{equation*}
Y_{n,0}=\sum_{j=1}^{n}\left( Y_{0},e_{j}\right)
_{V}e_{j}=\sum_{j=1}^{n}\left( Y_{0},\widetilde{e}_{j}\right) _{W}\widetilde{e}_{j},
\end{equation*}%
The Parseval's identity yields
\begin{equation*}
\left\Vert Y_{n}(0)\right\Vert _{V}\leq \left\Vert Y_{0}\right\Vert
_{V}\quad \text{ and }\quad \left\Vert Y_{n}(0)
\right\Vert _{W}\leq \left\Vert Y_{0}\right\Vert _{W}.  
\end{equation*}
The equation (\ref{y1}) can be written as a system of stochastic ordinary differential
equations in $\mathbb{R}^{n}$ with locally Lipschitz nonlinearities. From 
classical results there exists a local-in-time solution $\ Y_{n}$ that  is an adapted stochastic process with values in  $C([0,T_{n}],W_{n})$.

The existence of a   global-in-time solution
follows from the uniform estimates on $n=1,2,...,$ that will be deduced in the
next lemma \ (a similar reasoning can be found in \cite{ABW10}, \cite{CC18}, \cite{RS12}). 

\bigskip

\begin{lemma}
\label{existence_state} Let us assume \eqref{LG_1}-\eqref{R_1}. Then 
 the problem (\ref{y1}) admits a unique solution $Y_{n}\in L^{2}(\Omega
,L^{\infty }(0,T;W))$, which verifies the following estimates
\begin{eqnarray}
\mathbb{E}\sup_{s\in \lbrack 0, t]}\left\Vert
Y_{n}(s)\right\Vert _{V}^{2} &
+& 8\nu \mathbb{E}\int_{0}^{t}\left\Vert DY_{n}\right\Vert _{2}^{2}ds\,
+\frac{\beta}{4}\mathbb{E}\int_{0}^{
t}\|A_n\|_4^4\,ds  \notag \\
&\leq &C\left( 1+\mathbb{E}\left\Vert Y_{0}\right\Vert _{V}^{2}+\mathbb{E}%
\Vert U\Vert _{L^{2}(0,t;L^{2}(\mathcal{O}))}^{2}\right),  \quad   \forall t\in \lbrack 0,T], 
   \label{ineq1}
\end{eqnarray}
\begin{eqnarray}
\label{11:33_04:08} 
\mathbb{E}\|Y_n\|_{L^4(0,t;W^{1,4}(\mathcal{O}))}^4 \leq  C\left( 1+\mathbb{E}\left\Vert Y_{0}\right\Vert _{V}^{2}+\mathbb{E}%
\Vert U\Vert _{L^{2}(0,t;L^{2}(\mathcal{O}))}^{2}\right), \quad\forall t\in[0,T],
\end{eqnarray}
where $C$ is a positive  constant independent of $n$.
\end{lemma}

\textbf{Proof.} For each $n\in \mathbb{N}$, we define the following sequence of stopping times
\begin{equation*}
\tau _{M}^{n}=\inf \{t\geq 0:\Vert Y_{n}(t)\Vert_{V}\geq M\}\wedge
T_{n}, \quad M\in\mathbb{N}.
\end{equation*}
Let us set
\begin{equation}
\label{R_11}
f(Y_{n}):=\nu \Delta Y_{n}-
(Y_{n}\cdot \nabla)\upsilon_n-\sum_j\upsilon_n^j\nabla Y_{n}^j+(\alpha_1+\alpha_2)\,{\rm div}
\left(A_n^2\right)+\beta{\rm div}\left(|A_n|^2A_n\right)
+U.
\end{equation}
Using \eqref{ip_V}, and considering  in  (\ref{y1}) the  test functions $\phi =e_{i}$,  $i=1,\dots ,n$, we write
\begin{equation}
d\left( Y_{n},e_{i}\right) _{V}=\left( f(Y_{n}),e_{i}\right) \,dt+\left(
\sigma(t,Y_{n}),e_{i}\right) \,d\mathcal{W}_{t}.  \label{y3}
\end{equation}
Applying the It\^o formula, we deduce
\begin{equation*}
\label{24-05-1}
d\left( Y_{n},e_{i}\right) _{V}^{2}=2\left( Y_{n},e_{i}\right) _{V}\left(
f(Y_{n}),e_{i}\right) \,dt+2\left( Y_{n},e_{i}\right) _{V}\left(
\sigma(t,Y_{n}),e_{i}\right) \,d\mathcal{W}_t+|\left( \sigma\left( t,Y_{n}\right) ,e_{i}\right)|
^{2}\,dt.
\end{equation*}
 Summing 
over $i=1,\dots ,n,$ we derive
\begin{equation}
\label{06-17-38}
d\left\Vert Y_{n}\right\Vert _{V}^{2}=2\left( f(Y_{n}),Y_{n}\right)
\,dt+2\left( \sigma(t,Y_{n}),Y_{n}\right) \,d\mathcal{W}_t+\sum_{i=1}^{n}|\left(\sigma\left(
t,Y_{n}\right) ,e_{i}\right)|^{2}\,dt.
\end{equation}
We have
\begin{align}
\left( f(Y_{n}),Y_{n}\right) & =-2\nu \left\Vert DY_{n}\right\Vert
_{2}^{2}
-\bigl((Y_{n}\cdot \nabla)\upsilon_n+\sum_j\upsilon_n^j\nabla Y_{n}^j,Y_n\bigr)\notag \\
&+\left((\alpha_1+\alpha_2)\,{\rm div}
\left(A_n^2\right), Y_n\right)
+\left(\beta{\rm div}\left(|A_n|^2A_n\right), Y_n\right)
 +\left( U,Y_{n}\right)  \notag \\
& =I_1+I_2+I_3+I_4+I_5.
\label{equation_etat_temps_n1}
\end{align}
By the symmetry of the trilinear functional \eqref{tf_1}, we obtain
\begin{align}
\label{11-10-06}
I_2&=-\bigl((Y_{n}\cdot \nabla)\upsilon_n+\sum_j\upsilon_n^j\nabla Y_{n}^j,Y_n\bigr)
=-b(Y_{n},\upsilon_n,Y_{n})-b(Y_{n},Y_{n},\upsilon_n)\notag
\\&=-b(Y_{n},\upsilon_n,Y_{n})+b(Y_{n},\upsilon_n,Y_{n})=0.
\end{align}
Taking into account the boundary conditions  $Y_n=(Y_n\cdot \tau)\tau$,  $(n\cdot A_n)\cdot \tau=0$ on $\Gamma$ and the symmetry of $\nabla Y_n$, the divergence theorem gives
\begin{align}
\label{11-10-03}
I_4=\int_{\mathcal{O}}\beta \,{\rm div}\left(|A_n|^2A_n\right)\cdot Y_n
&=\beta\int_{\Gamma}
|A_n|^2(Y_n\cdot \tau)\,(n\cdot A_n)\cdot \tau-\beta\int_{\mathcal{O}} |A_n|^2A_n\cdot \nabla Y_n\notag\\
&=-\frac{\beta}{2}\|A_n\|_4^4.
\end{align}
Taking  $\epsilon=\frac{\beta}{4}$ in
\eqref{08-11-48}, we obtain 
 \begin{align}
 \label{11-10-04}
|I_3|&\leq  \frac{\beta}{4} \|A_n^2\|^2_2+\frac{(\alpha_1+\alpha_2)^2}{4\beta}\|A_n\|^2_2.
\end{align}
In addition, we have 
\begin{equation}
\label{11-10-05}
|I_5|=|(U,Y_n)|
\leq \frac{1}{2} \|Y_{n}\|_2^2+\frac{1}{2} \| U\|_2^2.
\end{equation}
Therefore, introducing \eqref{11-10-06}-\eqref{11-10-05} in 
\eqref{06-17-38}, we derive
\begin{align}
\label{02_12_61}
d\left\Vert Y_{n}\right\Vert _{V}^{2}&+\frac{\beta}{2}\|A_n\|_4^4\,dt+4\nu \left\Vert
DY_{n}\right\Vert _{2}^{2}\,dt
\leq \frac{(\alpha_1+\alpha_2)^2}{2\beta}\|A_n\|^2_2\,dt\notag\\
&+ (\| U\|_2^2+\|Y_{n}\|_2^2)\,dt+2\left( \sigma(t,Y_{n}),Y_{n}\right) \,d\mathcal{W}_t
+\sum_{i=1}^{n}|\left( \sigma\left(
t,Y_{n}\right) ,e_{i}\right)|^{2}\,dt.
\end{align}
We  write
\begin{align}
\label{10:07_04:08}
d\left\Vert Y_{n}\right\Vert _{V}^{2}&+\frac{\beta}{2}\|A_n\|_4^4\,dt+4\nu \left\Vert
DY_{n}\right\Vert _{2}^{2}\,dt
\leq C(\beta, \alpha_1,\alpha_2)\|Y_n\|^2_V\,dt\notag\\
&+ \| U\|_2^2\,dt+2\left( \sigma(t,Y_{n}),Y_{n}\right) \,d\mathcal{W}_t
+\sum_{i=1}^{n}|\left( \sigma\left(
t,Y_{n}\right) ,e_{i}\right)|^{2}\,dt.
\end{align}

Denoting by  $\widetilde{\sigma}_{n}$  the solution of the generalized Stokes problem (\ref{GS_NS}) for $
f=\sigma(t,Y_{n})$, we have 
\begin{equation*}
(\widetilde{\sigma}_{n},e_{i})_{V}=(\sigma(t,Y_{n}),e_{i})\quad \text{for }i=1,\dots ,n,\
\end{equation*}%
then   \eqref{LG_1}, \eqref{Korn_*} and 
Young's inequality give
\begin{align}
\sum_{i=1}^{n}|\left( \sigma\left( t,Y_{n}\right) ,e_{i}\right)| ^{2}
&=||\widetilde{\sigma}_{n}||_{V}^{2}\leq C||\sigma\left( t,Y_{n}\right) ||_{2}^{2}\leq
CL(1+\|Y_{n}\|)^\gamma_{W^{1,4}}\notag\\
&\leq CL+CL(K_*)^\gamma\|A_n\|^\gamma_4\leq C(L,\gamma,\beta)
+\frac{\beta}{4}\|A_n\|^4_4. 
 \label{est}
\end{align}
For any $t\in \lbrack 0,T]$, integrating the inequality 
 \eqref{10:07_04:08}
on  $(0,s)$,  $s\in[0,\tau _{M}^{n}\wedge t]$ and using 
\eqref{est}, we derive
\begin{align}
\label{R_3}
\left\Vert Y_{n}(s)\right\Vert _{V}^{2}&+\frac{\beta}{4}\int_0^s\|A_n\|_4^4\,dr+4\nu\int_0^s \left\Vert
DY_{n}\right\Vert _{2}^{2}\,dr
\leq \|Y_n(0)\|_V^2+C(\beta, \alpha_1,\alpha_2, T)\notag\\
&+C(\beta, \alpha_1,\alpha_2)\int_0^s\|Y_n\|^2_V\,dr+ \int_0^s\| U\|_2^2\,dr+2\int_0^t\left( \sigma(t,Y_{n}),Y_{n}\right) \,d\mathcal{W}_r.
\end{align}
On the other hand, the Burkholder-Davis-Gundy inequality, \eqref{LG_1},  \eqref{Korn_*} and  the Young inequality yield
\begin{align}
&\mathbb{E}\sup_{s\in \lbrack 0,\tau _{M}^{n}\wedge t]}\left\vert
\int_{0}^{s}\left( \sigma\left( r,Y_{n}\right) ,Y_{n}\right) \,dW_{r}\right\vert
 \leq C\mathbb{E}\left( \int_{0}^{\tau _{M}^{n}\wedge t}\left\vert \left(
\sigma\left( s,Y_{n}\right) ,Y_{n}\right) \right\vert ^{2}\,ds\right) ^{\frac{1}{2}}\notag \\
&
 \leq C\mathbb{E}\left( \int_{0}^{\tau _{M}^{n}\wedge t}
\| \sigma\left( s,Y_{n}\right)\|^2_2\|Y_{n}\|_2^{2}\,ds\right) ^{\frac{1}{2}}\leq C\mathbb{E}\left( \int_{0}^{\tau _{M}^{n}\wedge t}
L(1+\|Y_{n}\|^\gamma_{W^{1,4}})\|Y_{n}\|_V^{2}\,ds\right) ^{\frac{1}{2}}\notag\\
&
 \leq
 \frac{1}{4}\mathbb{E}\sup_{s\in \lbrack 0,\tau _{M}^{n}\wedge t]}\|Y_{n}\|_V^{2}+C^2LT+  \mathbb{E}\int_{0}^{\tau _{M}^{n}\wedge t}
C^2L(K_*)^\gamma\|A_{n}\|_4^{\gamma}\,ds\notag\\
&
 \leq
 \frac{1}{4}\mathbb{E}\sup_{s\in \lbrack 0,\tau _{M}^{n}\wedge t]}\|Y_{n}\|_V^{2}+C(L,T,K_*,\beta)+\frac{\beta}{16} \mathbb{E} \int_{0}^{\tau _{M}^{n}\wedge t}
\|A_{n}\|_4^{4}\,ds.
\label{R_2}
\end{align}
Taking the supremum on $s\in[0,\tau _{M}^{n}\wedge t]$ and the expectation
in \eqref{R_3} and incorporating the estimate \eqref{R_2}, we obtain
\begin{align}
&\frac{1}{2}\mathbb{E}\sup_{s\in \lbrack 0,\tau _{M}^{n}\wedge t]}\Vert
Y_{n}(s)\Vert _{V}^{2} +4\nu\mathbb{E}\int_{0}^{\tau _{M}^{n}\wedge t}
 \left\Vert DY_{n}\right\Vert _{2}^{2} \,ds
+\frac{\beta}{8}\mathbb{E}\int_{0}^{\tau _{M}^{n}\wedge t}\|A_n\|_4^4\,dr\notag\\
& \leq C(\beta, \alpha_1,\alpha_2,T)+\mathbb{E}\left\Vert
Y_{0}\right\Vert _{V}^{2} 
+\mathbb{E}\int_{0}^{ t}\left\Vert U\right\Vert _{2}^{2}\,ds
+ C(\beta, \alpha_1,\alpha_2)\mathbb{E%
}\int_{0}^{ t}\sup_{r\in \lbrack 0, \tau _{M}^{n}\wedge s]}\Vert
Y_{n}(r)\Vert _{V}^{2}\,ds.
\end{align}
 Then the function
\begin{equation*}
f(t)=\frac{1}{2}\mathbb{E}\sup_{s\in \lbrack 0,\tau _{M}^{n}\wedge t]}\Vert
Y_{n}(s)\Vert _{V}^{2} 
+4\nu\mathbb{E}\int_{0}^{\tau _{M}^{n}\wedge t}
 \left\Vert DY_{n}\right\Vert _{2}^{2} \,ds+\frac{\beta}{8}\mathbb{E}\int_{0}^{\tau _{M}^{n}\wedge
t}\|A_n\|_4^4\,ds
\end{equation*}%
fulfils the Gronwall's inequality
\begin{equation*}
f(t)\leq C +\mathbb{E}\left\Vert
Y_0\right\Vert _{V}^{2}+\mathbb{E}\int_{0}^{t}\left\Vert
U\right\Vert _{2}^{2}\,ds +C\int_{0}^{t}f(s)ds,
\end{equation*}%
which implies
\begin{eqnarray}
\mathbb{E}\sup_{s\in \lbrack 0,\tau _{M}^{n}\wedge t]}\left\Vert
Y_{n}(s)\right\Vert _{V}^{2} &+& 8\nu \mathbb{E}\int_{0}^{\tau _{M}^{n}\wedge
t}\left\Vert DY_{n}\right\Vert _{2}^{2}ds\,
+\frac{\beta}{4}\mathbb{E}\int_{0}^{\tau _{M}^{n}\wedge
t}\|A_n\|_4^4\,ds  \notag \\
&\leq &C\left( 1+\mathbb{E}\left\Vert Y_{0}\right\Vert _{V}^{2}+\mathbb{E}%
\Vert U\Vert _{L^{2}(0,t;L^{2}(\mathcal{O}))}^{2}\right) .  \label{IMP_15}
\end{eqnarray}
Then there exists a constant  $C$  independent of $M$ and $n$ such that
\begin{eqnarray}
\mathbb{E}\sup_{s\in \lbrack 0,\tau _{M}^{n}\wedge t]}\left\Vert
Y_{n}(s)\right\Vert _{V}^{2} \leq C, \quad \forall t\in [0,T].
\end{eqnarray}
 Let us fix $n\in \mathbb{%
N}$, writing
\begin{align}
\mathbb{E}\sup_{s\in \lbrack 0,\tau _{M}^{n}\wedge T]}\left\Vert
Y_{n}(s)\right\Vert _{V}^{2}& =\mathbb{E}\left( \sup_{s\in \lbrack
0,\tau _{M}^{n}\wedge T]}1_{\{\tau _{M}^{n}<T\}}\ \left\Vert
Y_{n}(s)\right\Vert _{V}^{2}\right)  +\mathbb{E}\left( \sup_{s\in \lbrack 0,\tau _{M}^{n}\wedge T]}1_{\{\tau
_{M}^{n}\geq T\}}\ \left\Vert Y_{n}(s)\right\Vert _{V}^{2}\right)  \notag
\\
& \geq \mathbb{E}\left( \max_{s\in \lbrack 0,\tau _{M}^{n}]}1_{\{\tau
_{M}^{n}<T\}}\ \left\Vert Y_{n}(s)\right\Vert _{V}^{2}\right)
\geq M^{2}P\left( \tau _{M}^{n}<T\right) ,  \label{fc}
\end{align}
we deduce that 
$
P\left( \tau _{M}^{n}<T\right)\leq \frac{C}{M^{2}}.
$
This means that $\tau _{M}^{n}\rightarrow T$ in
probability, as $M\rightarrow \infty $.
Then there exists a subsequence $%
\{\tau _{M_{k}}^{n}\}$ of $\{\tau _{M}^{n}\}$ (that may depend on $n$) such
that
\begin{equation*}
\tau _{M_{k}}^{n}\rightarrow T\quad \text{ a.e.  }\text{ as }k\rightarrow \infty .
\end{equation*}%
Since $\tau _{M_{k}}^{n}\leq T_{n}\leq T$, we deduce that $T_{n}=T$, so $%
Y_{n}$ is a global-in-time solution of the stochastic differential equation (%
\ref{y1}). In addition for fixed $n$, the monotonicity of the sequence $\left\{ \tau _{M}^{n}\right\} $
 allows to apply the  monotone convergence theorem in order to pass to the limit, as $M\rightarrow \infty $,  in the
inequality (\ref{IMP_15}) in order to obtain (\ref{ineq1}).

\begin{eqnarray}
\mathbb{E}\sup_{s\in[0, t]}\left\Vert
Y_{n}(s)\right\Vert _{V}^{2} &
+& 8\nu \mathbb{E}\int_{0}^{
t}\left\Vert DY_{n}\right\Vert _{2}^{2}ds\,
+\frac{\beta}{4}\mathbb{E}\int_{0}^{
t}\|A_n\|_4^4\,ds  \notag \\
&\leq &C\left( 1+\mathbb{E}\left\Vert Y_{0}\right\Vert _{V}^{2}+\mathbb{E}%
\Vert U\Vert _{L^{2}(0,t;L^{2}(\mathcal{O}))}^{2}\right) .  \label{IMP_1}
\end{eqnarray}
This   inequality  gives 
\begin{eqnarray}
\mathbb{E}\int_{0}^{
t}\|A_n\|_4^4\,ds &\leq C(\beta, \alpha_1, \alpha_2)\left( 1+\mathbb{E}\left\Vert Y_{0}\right\Vert _{V}^{2}+\mathbb{E}%
\Vert U\Vert _{L^{2}(0,t;L^{2}(\mathcal{O}))}^{2}\right), \quad 
\forall t\in[0,T],
\end{eqnarray}
that together with Lemma 2.1 yields
\begin{align} 
\mathbb{E}\|Y_n\|_{L^4(0,t;W^{1,4}(\mathcal{O}))}^4 &= \mathbb{E}  \int_0^t \|Y_n\|_{W^{1,4}}^4 ds 
\leq (K_*)^4 \mathbb{E} \int_{0}^{t}\|A_n\|_4^4\,ds\notag\\ 
& \leq  C(\beta, \alpha_1, \alpha_2)\left( 1+\mathbb{E}\left\Vert Y_{0}\right\Vert _{V}^{2}+\mathbb{E}%
\Vert U\Vert _{L^{2}(0,t;L^{2}(\mathcal{O}))}^{2}\right), \quad 
\forall t\in[0,T].
\end{align}
The H\"older's inequality also gives
\begin{eqnarray} 
\mathbb{E}\|Y_n\|_{L^4(0,t;W^{1,4}(\mathcal{O}))} \leq 
 C(\beta, \alpha_1, \alpha_2)\left( 1+\mathbb{E}\left\Vert Y_{0}\right\Vert _{V}^{2}+\mathbb{E}%
\Vert U\Vert _{L^{2}(0,t;L^{2}(\mathcal{O}))}^{2}\right)^{\frac{1}{4}}, \quad \forall t\in[0,T].
\end{eqnarray}

\bigskip\ $\hfill \hfill \blacksquare $

\begin{lemma}
\label{Ex_Int}
Assume  \eqref{LG_1}-\eqref{R_1}.  
Then we have 
\begin{align}
\label{17:10_04:08}
\mathbb{E}\text{e}^{\frac{\lambda \beta}{16(K_*)^4}\int_0^t\|Y_n\|_{W^{1,4}}^4 \,ds}
<C\,\mathbb{E}\text{e}^{\lambda\bigl(\int_0^T\|U\|_2^2ds+\|Y_0\|^2_V\bigr)}
, \quad \forall t\in[0,T],
\end{align}
where  $C$ is a positive constant independent of $n$, and  $K_*$  is defined 
by \ \eqref{Korn_*}.
\end{lemma}
\textbf{Proof.} Let us consider the inequality \eqref{R_3} and write
\begin{align}
\label{R_4}
\left\Vert Y_{n}(t)\right\Vert _{V}^{2}&+\frac{\beta}{4}\int_0^t\|A_n\|_4^4\,ds+4\nu\int_0^t \left\Vert
DY_{n}\right\Vert _{2}^{2}\,ds
\leq \|Y_n(0)\|_V^2+ \int_0^t\| U\|_2^2\,ds
\notag\\
&+C(\beta, \alpha_1,\alpha_2, T)+C(\beta, \alpha_1,\alpha_2)\int_0^t\|Y_n\|^2_V\,ds+2\int_0^t\left( \sigma(s,Y_{n}),Y_{n}\right) \,d\mathcal{W}_s.
\end{align}
Multiplying by  $\frac{\lambda}{2}$ and knowing that $W^{1,4}(\mathcal{O})\hookrightarrow H^1(\mathcal{O})$, we deduce
\begin{align*}
\frac{\lambda}{2}\left\Vert Y_{n}(t)\right\Vert _{V}^{2}&+\frac{\lambda\beta}{8}\int_0^t\|A_n\|_4^4\,ds+2\lambda\nu\int_0^t \left\Vert
DY_{n}\right\Vert _{2}^{2}\,ds
\leq \frac{\lambda}{2}\left(\|Y_n(0)\|_V^2+ \int_0^t\| U\|_2^2\,ds\right)
\notag\\
&+C(\beta, \alpha_1,\alpha_2, T, K_*)+C(\beta, \alpha_1,\alpha_2)
\int_0^t\|Y_n\|^2_{W^{1,4}}\,ds+\lambda\int_0^t\left( \sigma(s,Y_{n}),Y_{n}\right) \,d\mathcal{W}_s.
\end{align*}
The  Korn inequality \eqref{Korn_*} gives
$$
\frac{\lambda\beta}{8(K_*)^4}\|Y_n\|_{W^{1,4}}^4\leq \frac{\lambda\beta}{8}\|A_n\|_{4}^4;
$$
therefore we have
\begin{align}
\label{R_8}
\frac{\lambda\beta}{8(K_*)^4}\int_0^t\|Y_n\|_{W^{1,4}}^4\,ds
&\leq \frac{\lambda}{2}\left(\|Y_n(0)\|_V^2+ \int_0^t\| U\|_2^2\,ds\right)
+C(\beta, \alpha_1,\alpha_2, T, K_*)
\notag\\
&+C(\beta, \alpha_1,\alpha_2)
\int_0^t\|Y_n\|^2_{W^{1,4}}\,ds+\lambda\int_0^t\left( \sigma(s,Y_{n}),Y_{n}\right) \,d\mathcal{W}_s.
\end{align}
Let us notice that with the help of \eqref{LG_1}, the Sobolev embedding $W^{1,4}(\mathcal{O})\hookrightarrow H$ and the Young's inequality, for any $\delta>0$, we  can verify that
\begin{align*}
\int_0^t\lambda^2\left( \sigma(s,Y_{n},Y_{n}\right)^2 \,ds&\leq
\int_0^t\lambda^2\| \sigma(s,Y_{n}\|_2^2\|Y_{n}\|_2^2 \,ds
\leq
\int_0^t\lambda^2L(1+\|Y_{n})\|_{W^{1,4}}^\gamma)\|Y_{n}\|_2^2 \,ds
\notag\\
&\leq
\int_0^t\lambda^2L\|Y_{n}\|_{W^{1,4}}^2 \,ds
+\int_0^t\lambda^2L\|Y_{n})\|_{W^{1,4}}^{\gamma+2}\,ds\notag\\
&\leq C(\lambda,L,\delta,T)+\frac{\delta}{2}
\int_0^t\|Y_{n}\|_{W^{1,4}}^4 \,ds;
\end{align*}
which implies
$$
  -\frac{\delta}{2}
\int_0^t\|Y_{n}\|_{W^{1,4}}^4 \,ds-C(\lambda,L,\delta,T)\leq- \int_0^t\lambda^2\left( \sigma(s,Y_{n},Y_{n}\right)^2 \,ds.
$$
Adding this relation to \eqref{R_8}, we write
\begin{align}
\label{R_9}
\frac{\lambda\beta}{8(K_*)^4}\int_0^t\|Y_n\|_{W^{1,4}}^4\,ds&-\frac{\delta}{2}\int_0^t\|Y_n\|_{W^{1,4}}^4\,ds
\leq \frac{\lambda}{2}\left(\|Y_n(0)\|_V^2+ \int_0^t\| U\|_2^2\,ds\right)\notag\\
&
+C(\beta, \alpha_1,\alpha_2,\delta, T)
+C(\beta, \alpha_1,\alpha_2)
\int_0^t\|Y_n\|^2_{W^{1,4}}\,ds\notag\\
&+\lambda\int_0^t\left( \sigma(s,Y_{n}),Y_{n}\right) \,d\mathcal{W}_s- \int_0^t\lambda^2\left( \sigma(s,Y_{n},Y_{n}\right)^2 \,ds.
\end{align}
Once again, the Young inequality gives
$$C(\beta, \alpha_1,\alpha_2)
\int_0^t\|Y_n\|^2_{W^{1,4}}\,ds\leq C(\beta, \alpha_1,\alpha_2, \delta)+\frac{\delta}{2}
\int_0^t\|Y_n\|^4_{W^{1,4}}\,ds.$$
Introducing this estimate in \eqref{R_9} and next taking $\delta=\frac{\lambda\beta}{16(K_*)^4}$,
it follows that 
\begin{align}
\label{R_10}
\frac{\lambda\beta}{16(K_*)^4}\int_0^t\|Y_n\|_{W^{1,4}}^4\,ds&
\leq \frac{\lambda}{2}\left(\|Y_n(0)\|_V^2+ \int_0^t\| U\|_2^2\,ds\right)
+C(\beta, \alpha_1,\alpha_2, T)
\notag\\
&+\lambda\int_0^t\left( \sigma(s,Y_{n}),Y_{n}\right) \,d\mathcal{W}_s- \int_0^t\lambda^2\left( \sigma(s,Y_{n},Y_{n}\right)^2 \,ds.
\end{align}
Now, we take the exponential, the expectation and the H\"older inequality in order to deduce that
\begin{align*}
\mathbb{E}\text{e}^{\frac{\lambda\beta}{16(K_*)^4}\int_0^t\|Y_n\|_{W^{1,4}}^4\,ds }
&\leq C(\beta, \alpha_1,\alpha_2, T) \left(\mathbb{E}\text{e}^{\lambda \left(\|Y_0\|_V^2+ \int_0^t\| U\|_2^2\,ds\right)}\right)^{\frac{1}{2}}
\notag\\&
\qquad\left(\mathbb{E}\text{e}^{\int_0^t\left(2\lambda \sigma(s,Y_{n}),Y_{n}\right) \,d\mathcal{W}_s- \frac{1}{2}\int_0^t\left(2\lambda \sigma(s,Y_{n},Y_{n}\right)^2 \,ds}\right)^{\frac{1}{2}}
.
\end{align*}
Since the stochastic process inside the second expectation is a supermartingale its expectation is less or equal to $1$, hence we obtain \eqref{17:10_04:08}.

\bigskip\ $\hfill \hfill \blacksquare $

\begin{lemma}
\label{existence_state} Assume \eqref{LG_1}-\eqref{R_1}.
 Then  
the unique solution $Y_{n}$ of the problem (\ref{y1})  verifies the following uniform estimate
\begin{equation}
\mathbb{E}\sup_{s\in \lbrack 0,t]}\left\Vert Y_{n}(s)\right\Vert _{W}^{p}\leq C , \qquad   \forall t\in \lbrack 0,T],  \label{ineq2}
\end{equation}
where $C$ is a  positive constant independent of $n$. 
\end{lemma}

\textbf{Proof.} For each $n\in\mathbb{N}$, let us  consider the sequence of   stopping times
defined by
\begin{equation*}
\tau _{M}^{n}=\inf \{t\geq 0:\Vert Y_{n}(t)\Vert_{W}\geq M\}
, \quad M\in\mathbb{N}.
\end{equation*}

We introduce the solutions $\tilde{f}_{n}$ and $\tilde{\sigma}_{n}$ of (\ref{GS_NS}) for $%
f=f_n:=f(Y_{n})$ (as in \eqref{R_11}) and $f=\sigma_n:=\sigma(t,Y_{n})$, respectively. Then
\begin{equation}
(\tilde{f}_{n},e_{i})_{V}=(f_n,e_{i}),\qquad (\tilde{\sigma}%
_{n},e_{i})_{V}=(\sigma_n,e_{i}).  \label{gs1}
\end{equation}%
Therefore
\begin{equation*}
d( Y_{n},e_{i}) _{V}=( \tilde{f}_{n},e_{i})
_{V}\,dt+( \tilde{\sigma}_{n},e_{i}) _{V}\,d\mathcal{W}_t.
\end{equation*}%
Multiplying by $\lambda _{i}$ and using \eqref{mu}, we obtain
\begin{equation*}
d\left( Y_{n},e_{i}\right) _{W}=(\tilde{f}_{n},e_{i})_{%
W}\,dt+(\tilde{\sigma}_{n},e_{i})_{W}\,d\mathcal{W}_t.
\end{equation*}%
The It\^{o} formula gives
\begin{equation*}
d\left( Y_{n},e_{i}\right) _{W}^{2}=2\left( Y_{n},e_{i}\right) _{%
W}(\tilde{f}_{n},e_{i})_{W}\,dt+2\left(
Y_{n},e_{i}\right) _{W}(\tilde{\sigma}_{n},e_{i})_{W
}\,d\mathcal{W}_t+|(\tilde{\sigma}_{n},e_{i})_{W}|^{2}\,dt.
\end{equation*}%
Now, multiplying  by $\frac{1}{\lambda _{i}}$ and summing over $%
i=1,\dots ,n$, we derive
\begin{equation}
\label{02-20-20}
d\left\Vert Y_{n}\right\Vert _{W}^{2}=2(\tilde{f}_{n},Y_{n})_{%
W}\,dt+2(\tilde{\sigma}_{n},Y_{n})_{W}\,d\mathcal{W}_t+\sum_{i=1}^{n}%
\frac{1}{\lambda _{i}}|(\tilde{\sigma}_{n},e_{i})_{W}|^{2}\,dt,
\end{equation}
which is equivalent to
\begin{align}
\label{18:41_05:08}
d\left\Vert Y_{n}\right\Vert _{W}^{2}&=2\left[
(f_{n},Y_{n})
 +(f_{n},\mathbb{P}\upsilon(Y_{n}))
 \right]\,dt+\|\tilde{\sigma}_{n}\|_W^{2}\,dt+2\left[(\sigma_{n},Y_{n})+
(\sigma_{n},\mathbb{P}\upsilon(Y_{n}))
 \right]\,d\mathcal{W}_t.
\end{align}
Let us recall from \eqref{equation_etat_temps_n1}-\eqref{11-10-05}  that
\begin{align}
\label{18:20_05:08}
2(f_{n},Y_{n})\leq 
-4\nu \left\Vert
DY_{n}\right\Vert _{2}^{2} -\frac{\beta}{2}\|A_n\|_4^4
+ \frac{(\alpha_1+\alpha_2)^2}{2\beta}\|A_n\|^2_2+\|Y_n\|^2_2
+\| U\|_2^2.
\end{align}

On the other hand,
considering the Sobolev inequality
$$\|y\|_6\leq C_1\|y\|_{H^1}\quad \forall y\in H^1, 
$$
and  using the estimates \eqref{09-13-02}-\eqref{09-13-04}
as in  \cite{BI06}, page 373, for  $\epsilon={\rm min}\left\{
 \frac{1}{20(C_1)^2}, \frac{\alpha_1}{40(C_1)^2}, \frac{\beta\alpha_1}{9(3\beta+5)}\right\}$
we derive 
\begin{align}
\label{16:10_05:08}
2(f_{n},\mathbb{P}\upsilon(Y_{n}))&\leq -\frac{\beta}{2}\|A_n\|_4^4-\frac{\alpha_1\beta}{2}\||A_n||\nabla A_n|\|^2_2-\frac{\alpha_1\beta}{4}\|\nabla(|A_n|^2)\|^2_2\notag\\
&+C(\nu, \beta,\alpha_1,\delta)\|Y_n\|_W^2
+C(\beta,\alpha_1)\|Y_n\|_\infty
\|Y_n\|^2_{W}\notag\\
&
+C(\beta,\alpha_1)\|Y_n\|_{V}^2
\|Y_n\|^2_{W}
+2\delta\|Y_n\|^4_{W^{1,4}}+\|U\|^2_2.
\end{align}
The  Sobolev inequalities
$$\|y\|_\infty\leq C_2\|y\|_{W^{1,4}},\quad \|y\|_{V}\leq C_3\|y\|_{W^{1,4}},\quad \forall y\in W^{1,4}, 
$$
 \eqref{14:27_04:08} and the Young's inequality
 allow to verify that 
$$
C(\beta,\alpha_1)\left(\|Y_n\|_{\infty}+\|Y_n\|_{V}^2\right)\|Y_n\|^2_{W}
\leq C(\beta,\alpha_1, \delta)\|Y_n\|_{W}^2+2\delta \|y\|_{W^{1,4}}^4 \|Y_n\|_{W}^2,\quad \forall \delta>0.
$$
Therefore, we have
\begin{align}
\label{R_10}
2(f_{n},\mathbb{P}\upsilon(Y_{n}))&\leq -\frac{\beta}{2}\|A_n\|_4^4-\frac{\alpha_1\beta}{2}\||A_n||\nabla A_n|\|^2_2-\frac{\alpha_1\beta}{4}\|\nabla(|A_n|^2)\|^2_2\notag\\
&+C(\nu, \beta,\alpha_1,\delta)\|Y_n\|_W^2 +2\delta\|Y_n\|^4_{W^{1,4}}\|Y_n\|^2_{W}
+2\delta\|Y_n\|^4_{W^{1,4}}+\|U\|^2_2.
\end{align}
Now, we choose $\delta$ such that $2D_1:=4\delta\leq \frac{\lambda \beta}{16p(K_*)^4}$
and introduce the function  $$
\xi_1(t)=\text{e}^{-2D_1\int_0^t\|Y_n\|_{W^{1,4}}^4ds}.
$$
We apply the  It\^o formula to determine the differential of the product $\xi_1(t)\|Y_n(t)\|_W^2$, namely  from the  equation \eqref{18:41_05:08} we derive 
\begin{align}
\label{R_14}
d\left(\xi_1(t)\left\Vert Y_{n}\right\Vert _{W}^{2}\right)
&=\xi_1(t)\left[
2(f_{n},Y_{n})
 +2(f_{n},\mathbb{P}\upsilon(Y_{n}))
 \right]\,dt+\xi_1(t)\|\tilde{\sigma}_{n}\|_W^{2}\,dt\notag\\
 &+\xi_1(t)\left[2(\sigma_{n},Y_{n})+
2(\sigma_{n},\mathbb{P}\upsilon(Y_{n}))
 \right]\,d\mathcal{W}_t\notag\\
 &-2D_1\xi_1(t)\left\Vert Y_{n}\right\Vert _{W}^{2}
 \|Y_n\|_{W^{1,4}}^4\,dt.
\end{align}
Using the It\^o formula once again for the function $\theta(x)=x^p$, 
and integrating   on $[0,s]$, $s\leq t\wedge \tau^n_M$, $t\in[0,T]$, 
we deduce
\begin{align}
\label{18:43_05:08}
\left(\xi_1(s)\left\Vert Y_{n}\right\Vert _{W}^{2}\right)^p&=
\| Y_{n}(0)\|_{W}^{2p}
+p\int_0^s\left(\xi_1(r)\left\Vert Y_{n}\right\Vert _{W}^{2}\right)^{p-1}\xi_1(r)\left[
2(f_{n},Y_{n})
 +2(f_{n},\mathbb{P}\upsilon(Y_{n}))
 \right]\,dr\notag\\
 &+p\int_0^s\left(\xi_1(r)\left\Vert Y_{n}\right\Vert _{W}^{2}\right)^{p-1}\xi_1(r)
 \|\tilde{\sigma}_{n}\|_W^{2}\,dr\notag\\
 &+p\int_0^s\left(\xi_1(r)\left\Vert Y_{n}\right\Vert _{W}^{2}\right)^{p-1}\xi_1(r)\left[2(\sigma_{n},Y_{n})+
2(\sigma_{n},\mathbb{P}\upsilon(Y_{n}))
 \right]\,d\mathcal{W}_r\notag\\
 &-2D_1p\int_0^s\left(\xi_1(r)\left\Vert Y_{n}\right\Vert _{W}^{2}\right)^{p-1}\xi_1(r)\left\Vert Y_{n}\right\Vert _{W}^{2}
 \|Y_n\|_{W^{1,4}}^4\,dr\notag\\
 &+2p(p-1)\int_0^s\left(\xi_1(r)\left\Vert Y_{n}\right\Vert _{W}^{2}\right)^{p-2}(\xi_1(r))^2\left[(\sigma_{n},Y_{n})+
(\sigma_{n},\mathbb{P}\upsilon(Y_{n}))
 \right]^2\,dr.
\end{align}
Next, using \eqref{18:20_05:08} and \eqref{R_10}
to estimate the right hand side, we obtain
\begin{align}
\label{18:44_05:08}
\left(\xi_1(s)\left\Vert Y_{n}\right\Vert _{W}^{2}\right)^p&\leq 
\| Y_{n}(0)\|_{W}^{2p}
+p\int_0^s\left(\xi_1(r)\left\Vert Y_{n}\right\Vert _{W}^{2}\right)^{p-1}\xi_1(r)
\biggl[
 C(\nu,\beta,\alpha_1, \delta)\|Y_n\|_W^2
+ 2\| U\|_2^2 \biggr]\,dr\notag\\
& + p\int_0^s\left(\xi_1(r)\left\Vert Y_{n}\right\Vert _{W}^{2}\right)^{p-1}\xi_1(r)
\left[D_1\|Y_n\|^4_{W^{1,4}}
\|Y_n\|^2_{W}+D_1\|Y_n\|^4_{W^{1,4}}
 \right]\,dr\notag\\
 &+p\int_0^s\left(\xi_1(r)\left\Vert Y_{n}\right\Vert _{W}^{2}\right)^{p-1}\xi_1(r)
 \|\tilde{\sigma}_{n}\|_W^{2}\,ds\notag\\
 &+2p\int_0^s\left(\xi_1(r)\left\Vert Y_{n}\right\Vert _{W}^{2}\right)^{p-1}\xi_1(r)\left[(\sigma_{n},Y_{n})+
(\sigma_{n},\mathbb{P}\upsilon(Y_{n}))
 \right]\,d\mathcal{W}_r\notag\\
 &-2D_1p\int_0^s\left(\xi_1(r)\left\Vert Y_{n}\right\Vert _{W}^{2}\right)^{p-1}\xi_1(r)\left\Vert Y_{n}\right\Vert _{W}^{2}\|Y_n\|_{W^{1,4}}^4
 \,dr\notag\\
 &+2p(p-1)\int_0^s\left(\xi_1(r)\left\Vert Y_{n}\right\Vert _{W}^{2}\right)^{p-2}(\xi_1(r))^2\left[(\sigma_{n},Y_{n})+
(\sigma_{n},\mathbb{P}\upsilon(Y_{n}))
 \right]^2\,dr.
\end{align}
Since $\| Y_{n}\|_{W}^{2p-2}\leq 1+\| Y_{n}\|_{W}^{2p}$, we deduce
\begin{align}
\label{18:45_05:08}
\left(\xi_1(s)\left\Vert Y_{n}\right\Vert _{W}^{2}\right)^p&\leq 
\| Y_{n}(0)\|_{W}^{2p}
+pC(\nu,\beta,\alpha_1, \delta)\int_0^s\left(\xi_1(r)
\left\Vert Y_{n}\right\Vert _{W}^{2}\right)^{p}
 \,dr\notag\\
& + pD_1\int_0^s\left(\xi_1(r)\right)^p\|Y_n\|^4_{W^{1,4}} \,dr\notag\\ 
 &
+ 2p\int_0^s\left(\xi_1(r)\left\Vert Y_{n}\right\Vert _{W}^{2}\right)^{p-1}\xi_1(r)\|U\|^2_2
\,dr\notag\\ 
 &+p\int_0^s\left(\xi_1(r)\left\Vert Y_{n}\right\Vert _{W}^{2}\right)^{p-1}\xi_1(r)
 \|\tilde{\sigma}_{n}\|_W^{2}\,dr\notag\\
 &+2p\int_0^s\left(\xi_1(r)\left\Vert Y_{n}\right\Vert _{W}^{2}\right)^{p-1}\xi_1(r)\left[(\sigma_{n},Y_{n})+
(\sigma_{n},\mathbb{P}\upsilon(Y_{n}))
 \right]\,d\mathcal{W}_s\notag\\
  &+2p(p-1)\int_0^s\left(\xi_1(r)\left\Vert Y_{n}\right\Vert _{W}^{2}\right)^{p-2}(\xi_1(r))^2\left[(\sigma_{n},Y_{n})+
(\sigma_{n},\mathbb{P}\upsilon(Y_{n}))
 \right]^2\,dr.
\end{align}
Taking into account that $ \|\tilde{\sigma}_{n}\|_W^{2}\leq C \|\sigma\|_2^2$, using   \eqref{LG_1},  the  embedding
 $W\hookrightarrow W^{1,4}(\mathcal{O})$ and the Young's inequality, we infer that
\begin{align*}
&p\int_0^s\left(\xi_1(r)\left\Vert Y_{n}\right\Vert _{W}^{2}\right)^{p-1}\xi_1(r)
 \|\tilde{\sigma}_{n}\|_W^{2}\,dr\notag\\
   &+2p(p-1)\int_0^s\left(\xi_1(r)\left\Vert Y_{n}\right\Vert _{W}^{2}\right)^{p-2}(\xi_1(r))^2\left[(\sigma_{n},Y_{n})+
(\sigma_{n},\mathbb{P}\upsilon(Y_{n}))
 \right]^2\,dr\notag\\
 &\leq C(p,T)+C(p)\int_0^s\left(\xi_1(r)
\left\Vert Y_{n}\right\Vert _{W}^{2}\right)^{p}
 \,dr.
\end{align*}
On the other hand, the Young's inequality \eqref{R_15} with
$r=\frac{p}{p-1}$ and $0\leq \xi_1(t)\leq 1$ give
$$
 2p\int_0^s\left(\xi_1(r)\left\Vert Y_{n}\right\Vert _{W}^{2}\right)^{p-1}\xi_1(r)\|U\|^2_2
\,dr
\leq 2(p-1) \int_0^s\left(\xi_1(r)\|Y_n\|_W^2\right)^{p}dr+2\int_0^s
  \| U\|_2^{2p}dr.
$$
Introducing these estimates in \eqref{18:45_05:08}, we deduce
\begin{align}
\label{R_17}
\left(\xi_1(s)\|Y_n(s)\|_W^2\right)^p
&\leq 
\|Y_0\|_W^{2p}
+C\biggl(\int_0^s\left(\xi_1(r)
\left\Vert Y_{n}\right\Vert _{W}^{2}\right)^{p}
 \,dr+\int_0^s\|Y_n\|^4_{W^{1,4}}\,dr
+\int_0^s
  \| U\|_2^{2p}
dr+1\biggr)
\notag\\
& +2p\int_0^s \left(\xi_1(r)\|Y_n\|_W^2\right)^{p-1}
\xi_1(r) \left[( \sigma(r,Y_{n}),\mathbb{P}\upsilon(Y_{n})+Y_{n})\,d\mathcal{W}_r\right].
\end{align}
The   Burkholder-Davis-Gundy inequality,   \eqref{LG_1} and the Young's inequality yield
\begin{align}
\mathbb{E}\sup_{s\in \lbrack 0,\tau^n _{M}\wedge t]}&
\left|
\int_0^s \left(\xi_1(r)\|Y_n\|_W^2\right)^{p-1}
\xi_1(r) \left(( \sigma(r,Y_{n}),\mathbb{P}\upsilon(Y_{n})+Y_{n})\,d\mathcal{W}_r\right)\right|
  \notag \\
& \leq C\mathbb{E}\left( 
\int_{0}^{\tau^n _{M}\wedge t}
\left(\left(\xi_1(s)\|Y_n\|_W^2\right)^{p-1}
\xi_1(s)\right)^{2}\| \sigma(s,Y_{n})\|^2_2\|Y_n\|_W^2
\,ds
\right) ^{\frac{1}{2}}
\notag \\
& \leq C\sqrt{L}\mathbb{E}\left( 
\int_{0}^{\tau^n _{M}\wedge t}
\left(\xi_1(s)\right)^{2p}\|Y_n\|_W^{4p-2}
+\int_{0}^{\tau^n _{M}\wedge t}
\left(\xi_1(s)\right)^{2p}\|Y_n\|_W^{4p-2+\gamma}
\right) ^{\frac{1}{2}}
\notag \\
& \leq  C\sqrt{2LT}+C\sqrt{2L}\mathbb{E}\left(
\int_{0}^{\tau^n _{M}\wedge t}\left(\xi_1(s)\|Y_n\|_W^2\right)^{2p}
\,ds
\right) ^{\frac{1}{2}}
\notag \\
& \leq   C\sqrt{2LT}+\frac{\eta}{2p}\mathbb{E}\left(\sup_{s\in \lbrack 0,\tau^n _{M}\wedge t]}
\left(\xi_1(s)\|Y_n\|_W^2\right)^{p}\right)\notag\\
&+C(L,\eta,p)\mathbb{E}
\int_{0}^{\tau^n _{M}\wedge t}\left(\xi_1(s)\|Y_n\|_W^2\right)^{p}
\,ds,
\label{s1}
\end{align}
for any $\eta>0.$ Here we  take $\eta=\frac{1}{2}.$ 
Considering the supremum on $s\in[0,\tau^n_{M}\wedge t]$ and the expectation in \eqref{R_17}, with the help of \eqref{s1} we derive  the following Gronwall's inequality 
\begin{align*}
&\frac{1}{2}\mathbb{E}\sup_{s\in \lbrack 0,\tau^n_{M}\wedge t]}\left(\xi_1(s)\|Y_n(s)\|_W^2\right)^p
\leq 
\|Y_0\|_W^{2p}
+C\biggl(
\int_0^t\mathbb{E}\sup_{r\in \lbrack 0,\tau^n_{M}\wedge s]} \left(\xi_1(r)\|Y_n(r)\|_W^2\right)^{p}\, ds
\notag\\
&
+\mathbb{E}\int_0^{\tau^n_{M}\wedge t} \|Y_n\|^4_{W^{1,4}}\,dr+\mathbb{E}\int_0^{\tau^n_{M}\wedge t}
 \| U\|_2^{2p}\,
dr+1\biggr).
\end{align*}
Therefore, we obtain
\begin{align}
\label{12-19-44}
\mathbb{E}\sup_{s\in[0,\tau^n_M\wedge t]}\left(\xi_1(s)\|Y_n(s)\|_W^2\right)^p
&\leq C\left(1+ \mathbb{E}\int_0^{t}
\| U\|_2^{2p}\,dr+
  \mathbb{E}\int_0^{ t}\|Y_n\|^4_{W^{1,4}}dr
   \right).
\end{align}
The estimates \eqref{11:33_04:08} and \eqref{R_1}   yield
\begin{align*}
\mathbb{E}\sup_{s\in[0,\tau^n_M \wedge t]}\left(\xi_1(s)\|Y_n(s)\|_W^2\right)^p
&\leq  C
\end{align*}
with $C$ independent of $n$ and $M$. We verify that 
for  $n$ fixed, $\tau
_{M}^{n}\rightarrow T$ in probability, as $M\rightarrow \infty $. Then,
there exists a subsequence $\{\tau _{M_{k}}^{n}\}$ of $\{\tau _{M}^{n}\}$
(that may depend on $n$) such that $\tau _{M_{k}}^{n}\rightarrow T$ \ for a.
e. $\omega \in \Omega $, as $k\rightarrow \infty $. Using the monotone
convergence theorem, we pass to the limit in (\ref{12-19-44}) as $k\rightarrow
\infty $, deriving the estimate 
\begin{align*}
\mathbb{E}\sup_{s\in[0,t]}\left(\xi_1(s)\|Y_n(s)\|_W^2\right)^p
&\leq C.
\end{align*}
The H\"older inequality gives 
\begin{align*}
\mathbb{E}\sup_{s\in \lbrack 0,t]}\left\Vert
Y_{n}(s)\right\Vert _{W}^{p}&\leq  \mathbb{E}\biggl[\biggl(\sup_{s\in \lbrack 0,t]}
(\xi_1(s))^{\frac{p}{2}}
\left\Vert
Y_{n}(s)\right\Vert _{W}^{p}\biggr)
(\xi_1(t))^{-\frac{p}{2}}\biggr]\notag\\
&\leq  \left(\mathbb{E}\sup_{s\in \lbrack 0,t]}
\left(\xi_1(s)\|Y_n(s)\|_W^2\right)^p\right)^{\frac{1}{2}}
\biggl(\mathbb{E}(\xi_1(t))^{-p}\biggr)^{\frac{1}{2}}\notag\\
&\leq \sqrt{C}
\biggl(\mathbb{E}
\text{e}^{2pD_1\int_0^t\|Y_n\|_{W^{1,4}}^4ds}
\biggr)^{\frac{1}{2}}.
\end{align*}
Using Lemma \ref{Ex_Int},  we deduce \eqref{ineq2}.

\bigskip

\subsection{Proof of Theorem \ref{the_1}.} 

In order to show the existence of the solution to the system $\eqref{equation_etat}$ it is convenient to write  the equation  $\eqref{equation_etat}_1$
in the following  form (see \cite{BI04}, page 3)
\begin{equation}
\label{eq_1}
\begin{array}{ll}
d(\upsilon (Y) )=
\bigl(-\nabla p+ \nu \Delta Y-
(Y\cdot \nabla)Y
+{\rm div}
N(Y)+{\rm div}S(Y)+U\bigr)dt+ \sigma(t,Y) \,d\mathcal{W}_{t},
 \end{array}
\end{equation}
with  the operators $S$ and $N$ defined in \eqref{operator1}-\eqref{operator2}.
The corresponding finite dimensional approximation  reads
\begin{eqnarray}
\label{11Ab1}
d(\upsilon (Y_n))=\left(-\nabla p_n+\nu \Delta Y_n-(Y_n\cdot\nabla) Y_n
+{\rm div}\,N(Y_n)+{\rm div}\,S(Y_n)+U\right)dt
+\sigma(t,Y_n)\,d\mathcal{W}_{t}
\end{eqnarray}

The proof of  Theorem \ref{the_1} is splitted into five steps. \vspace{2mm}\newline
\textit{Step 1. Convergences related with the projection operator.}
Let $P_{n}:W\rightarrow W_{n}$ be the orthogonal projection
defined by%
\begin{equation*}
P_{n}y=\sum_{j=1}^{n}\hat{c}_{j}\hat{e}_{j}\quad \text{ with }%
\quad \text{ }\hat{c}_{j}=\left( y,\hat{e}_{j}\right) _{%
W},\qquad \forall y\in W,
\end{equation*}%
where $\{\hat{e}_{j}=\frac{1}{\sqrt{\lambda _{j}}}e_{j}\}_{j=1}^{%
\infty }$ is the orthonormal basis of $\ W.$ \ \bigskip It is
easy to check that
\begin{equation*}
P_{n}y=\sum_{j=1}^{n}c_{j}e_{j}\quad \text{ with }\quad \text{ }c_{j}=\left(
y,e_{j}\right) _{V},\qquad \forall y\in W.
\end{equation*}

By Parseval's identity we have that
\begin{equation*}
||P_{n}y||_{V}\leq ||y||_{V},\text{ }\quad \forall y\in V,
\end{equation*}
\begin{equation*}
||P_{n}y||_{W}\leq ||y||_{W}\qquad \mbox{ and}
\qquad P_{n}y\longrightarrow y\qquad \mbox{ strongly in
} W,\quad \forall y\in W.
\end{equation*}

Considering an arbitrary $Z\in L^{q}(\Omega \times (0,T);W)),$
we have%
\begin{equation*}
||P_{n}Z||_{W}\leq ||Z||_{W}\qquad \mbox{ and
}\qquad P_{n}Z(\omega ,t)\rightarrow Z(\omega ,t)\qquad \mbox{ strongly in
}\ W,
\end{equation*}%
which are valid for $P$-a.e. $\omega \in \Omega $ and a.e. $t\in (0,T).$
Hence Lebesgue's dominated convergence theorem and the inequality
\begin{equation*}
||Z||_{V}\leq C||Z||_{W}\qquad \text{for any }Z\in W
\end{equation*}%
imply%
\begin{eqnarray}
P_{n}Z &\longrightarrow &Z\qquad \mbox{ strongly in
}\ L^{q}(\Omega \times (0,T),W)),  \notag \\
P_{n}Z &\longrightarrow &Z\qquad \mbox{ strongly in
}\ L^{q}(\Omega \times (0,T),V)).  \label{c02}
\end{eqnarray}%
\bigskip \textit{Step 2.} \textit{Passing to the limit in the weak sense.}
From Lemma \ref{existence_state}, we have
\begin{equation}
 \mathbb{E}\sup_{t\in \lbrack 0,T]}\left\Vert
Y_{n}(t)\right\Vert _{W}^{q}\leq C.  \label{as}
\end{equation}
  Then there exists 
 a subsequence of $Y_{n}$, still denoted by $Y_{n}$
such that
\begin{eqnarray}
Y_{n} &\rightharpoonup &Y\qquad \mbox{ *-weakly in
}\ L^{q}(\Omega ,L^{\infty }(0,T;W)).  \label{c1}
\end{eqnarray}%
Moreover, we have
\begin{eqnarray}
P_{n}Y &\longrightarrow &Y\qquad \mbox{ strongly in}
\ L^{q}(\Omega \times (0,T),W).  \label{c02Y}
\end{eqnarray}
Let us notice 
\begin{equation*}
|( S(y) ,\phi)| \leq C \|y\|_W^3\|\phi\|_2 \qquad \mbox{ for any
}\ y\in W\qquad \mbox{ and
}\ \phi\in H,
\end{equation*}
which implies that  $S:W\rightarrow H^{\ast }$ and 
$$
\|S(y)\|_{ H^{\ast }}\leq C \|y\|_W^3, \qquad \forall y\in W.
$$ 
Therefore 
\begin{equation}
\label{13-08-53}
\|S(Y_n)\|^2_{L^2(\Omega, L^2(0,T;H^{\ast }))}=\mathbb{E}\int_0^T\|S(Y_n)\|^2_{ H^{\ast }}\leq 
C\mathbb{E}\sup_{{t\in[0,T]}}\|Y_n(t)\|_W^6<C.
\end{equation}
We also have 
\begin{equation*}
|{\rm div\,}( S(y) ,\phi)| \leq C \|y\|_W^3\|\phi\|_V \qquad \mbox{ for any
}\ y\in W\qquad \mbox{ and
}\ \phi\in V,
\end{equation*}
then 
$$
\|{\rm div\,}S(y)\|_{ V^{\ast }}\leq C\|y\|_W^3, \qquad \forall y\in W,
$$
and 
\begin{equation}
\label{13-08-54}
\|{\rm div\,}S(y)\|^2_{L^2(\Omega, L^2(0,T;W^{\ast }))}\leq \|{\rm div\,}S(y)\|^2_{L^2(\Omega, L^2(0,T;V^{\ast }))}<C.
\end{equation}
The operator $N$ verifies
\begin{equation*}
|( N(y) ,\phi)| \leq C
 \|y\|_W^2\|\phi\|_2 \qquad \mbox{ for any
}\ y,\phi \in W\qquad \mbox{ and
}\ \phi\in V,
\end{equation*}
In addition
\begin{equation*}
|( {\rm div}\,N(y) ,\phi)| \leq C
 \|y\|_W^2\|\phi\|_W \qquad \mbox{ for any
}\ y,\phi \in W\qquad \mbox{ and
}\ \phi\in V,
\end{equation*}
which imply 
\begin{equation}
\label{13-08-55}
\|N(y)\|^2_{L^2(\Omega, L^2(0,T;H^{\ast }))}<C,
\end{equation}and 
\begin{equation}
\label{13-09-21}
\|{\rm div}\,N(y)\|^2_{L^2(\Omega, L^2(0,T;W^{\ast }))}\leq \|{\rm div}\,N(y)\|^2_{L^2(\Omega, L^2(0,T;V^{\ast }))}<C.
\end{equation}
Let us introduce the operator  $B$, defined by
$$
B(y):=-(y\cdot \nabla)y.
$$
We have 
\begin{equation}
|(B(y),\phi)|\leq C\|y\|_V^2\|\phi\|_V,
\label{14:31_07:08}
\end{equation}
then
\begin{equation}
\left\Vert B(Y_n)\right\Vert^2 _{L^{2}(\Omega, L^2(0,T;V^{\ast
}))}\leq C_{1}\left\Vert Y_n\right\Vert _{L^{4}(\Omega ,L^{\infty}(0,T;V))}^{2}<C.
\label{l}
\end{equation}
The diffusion operator is bounded.
Then there exist  operators $N^{\ast }(t)$, $S^{\ast }(t)$, $B^{\ast }(t)$, $\sigma^{\ast }(t)$ and a subsequence on $(n)$, that we still denote by $(n)$, such that as $n\to \infty$ we have 
\begin{eqnarray}
B(Y_{n}) &\rightharpoonup &B^{\ast }(t)\qquad \mbox{ weakly in
}\ L^{2}(\Omega \times (0,T),V^{\ast }),  \notag \\
N(Y_{n}) &\rightharpoonup &N^{\ast }(t)\qquad \mbox{ weakly in
}\ L^{2}(\Omega \times (0,T),H^{\ast }),  \notag \\
{\rm div}\,N(Y_{n}) &\rightharpoonup &{\rm div}\,N^{\ast }(t)\qquad \mbox{ weakly in
}\ L^{2}(\Omega \times (0,T),V^{\ast }),  \notag \\
S(Y_{n}) &\rightharpoonup &S^{\ast }(t)\qquad \mbox{ weakly in
}\ L^{2}(\Omega \times (0,T),H^{\ast }),  \notag \\
{\rm div}\,S(Y_{n}) &\rightharpoonup &{\rm div}\,S^{\ast }(t)\qquad \mbox{ weakly in
}\ L^{2}(\Omega \times (0,T),V^{\ast }),  \notag \\
\sigma(t,Y_{n}) &\rightharpoonup &\sigma^{\ast }(t)\qquad \mbox{ weakly in
}\ L^{2}(\Omega \times (0,T),(L^2(O))^m).  \label{c01}
\end{eqnarray}%
Therefore, passing to the limit with respect to the weak topology,
as $n\rightarrow \infty $, all terms in the  equation (\ref{y1}), we
derive that the limit function $Y$ satisfies \ the stochastic differential
equation
\begin{align}
d\left( \upsilon \left( Y\right) ,\phi \right) &=\left[ \left( \nu \Delta
Y+U,\phi \right) +\langle B^{\ast }(t),\phi \rangle
+\langle{\rm div\,}N^{\ast }(t),\phi \rangle 
+\langle {\rm div\,}S^{\ast }(t),\phi \rangle
\right] \,dt\vspace{2mm}%
+\left( \sigma^{\ast }(t),\phi \right) \,d\mathcal{W}_t,\notag\\
&\qquad \qquad \forall \phi \in
V.  \label{y12}
\end{align}

\textit{Step 3.} \textit{Passing to the limit in the strong  sense up to a stopping time.}
Let us introduce the following convenient sequence $(\tau _{M})$, $M\in \mathbb{N}$, of 
stopping times
\begin{equation*}
\tau _{M}=\inf \{t\geq 0:\Vert Y(t)\Vert _{W}
\geq M\}\wedge T.
\end{equation*}
\begin{proposition}
Let $Y_n$ be the solution of \eqref{11Ab1} and $P_{n}Y$ the orthogonal projection of the weak limit $Y$ on the space $W_n.$ Then for $M$ fixed we have
\begin{align}
\label{15:47_07:08}
& \mathbb{E} \left( \xi_2 (t\wedge \tau_M)||P_{n}Y(t\wedge \tau_M)-Y_{n}(t\wedge \tau_M)||_{V}^{2}\right)
+4\nu\mathbb{E}\int_0^{t\wedge \tau_M}\xi_2 (s) 
||D(P_{n}Y-Y_{n})||_{2}^{2} ds \notag\\
&+ 
\frac{\beta}{2}\mathbb{E}\int_0^{t\wedge \tau_M}\xi_2 (s)\int_{\mathcal{O}}(|A_n|^2-|A|^2)^2ds
+ \frac{\beta}{4}\mathbb{E}\int_0^{t\wedge \tau_M} \xi_2 (s) \int_{\mathcal{O}}(|A_n|^2+|A|^2)|A(Y_n-Y)|^2 ds\notag \\
&+\mathbb{E}\int_0^{t\wedge \tau_M} \xi_2 (s)\|P_n\tilde{\sigma}-P_n\tilde{\sigma^*}\|^2_Vds \to 0, \qquad\text{as}\quad n\to\infty,
\end{align}
where 
\begin{equation*}
\xi_2 (t)=e^{-D_{3}t-2D_{4}\int_{0}^{t}\left\Vert Y\right\Vert _{W
}ds}
\end{equation*}
and $D_{3}$, $D_{4}$ are specific constants to be defined later on.
\end{proposition}
\textbf{Proof.} 
Taking the difference between the equations (\ref{y1})  and (\ref{y12}), we write
\begin{eqnarray}
d\left( Y_{n}-P_{n}Y,e_{i}\right) _{V} &=&\left[ \left( \nu \Delta
(Y_{n}-Y),e_{i}\right) +\langle B(Y_{n})-B^{\ast }(t),e_{i}\rangle %
\right] \,dt \notag\\
&&+\left[ 
\langle {\rm div}\,N(Y_{n})-{\rm div}\,N^{\ast }(t),e_{i}\rangle
+
\langle {\rm div}\,S(Y_{n})-{\rm div}\,S^{\ast }(t),e_{i}\rangle
\right] \,dt \notag \\
&&\vspace{2mm}+\left( \sigma(t,Y_{n})-\sigma^{\ast }(t),e_{i}\right) \,d\mathcal{W}_t,
\label{diff}
\end{eqnarray}%
which holds for any $e_{i}\in W_{n},$ $i=1,...,n.$

The  It\^{o}'s formula gives
\begin{align*}
d(Y_{n}-P_{n}Y,e_{i})_{V}^{2}&=2\left( Y_{n}-P_{n}Y,e_{i}\right) _{V}\left[
\left( \nu \Delta (Y_{n}-Y),e_{i}\right) +\langle B(Y_{n})-B^{\ast
}(t),e_{i}\rangle \right] \,dt
\vspace{2mm} \\
&+2\left( Y_{n}-P_{n}Y,e_{i}\right) _{V}
\left[ 
\langle {\rm div}\,N(Y_{n})-{\rm div}\,N^{\ast }(t),e_{i}\rangle
+
\langle {\rm div}\,S(Y_{n})-{\rm div}\,S^{\ast }(t),e_{i}\rangle
\right] \,dt 
\vspace{2mm} \\
&+2\left( Y_{n}-P_{n}Y,e_{i}\right) _{V}\left( \sigma(t,Y_{n})-G^{\ast
}(t),e_{i}\right) \,d\mathcal{W}_t+|\left( \sigma(t,Y_{n})-\sigma^{\ast }(t),e_{i}\right)|
^{2}\,dt.
\end{align*}
Summing on  $i=1,\dots,n$, we obtain 
\begin{align}
&d\left( ||Y_{n}-P_{n}Y||_{V}^{2}\right) + 4\nu
||D(Y_{n}-P_{n}Y)||_{2}^{2} dt  \notag \\
& =
2\nu\left( \Delta(P_{n}Y-Y), Y_{n}-P_{n}Y\right)  \,dt
+
2\langle B(Y_{n})-B^{\ast }(t),Y_{n}-P_{n}Y\rangle \,dt
 \notag\\
&+2
\left[ 
\langle {\rm div}\left(N(Y_{n})-N^{\ast }(t)\right), Y_{n}-P_{n}Y\rangle
+
\langle {\rm div}\left(S(Y_{n})-S^{\ast }(t)\right), Y_{n}-P_{n}Y\rangle
\right] \,dt 
  \notag \\
& +\sum_{i=1}^{n}|\left( \sigma(t,Y_{n})-\sigma^{\ast }(t),e_{i}\right)|
^{2}\,dt+2\left( \sigma(t,Y_{n})-\sigma^{\ast }(t),Y_{n}-P_{n}Y\right) \,d\mathcal{W}_t.
\label{y13}
\end{align} 
Now, we write each term in the right hand side of this equation in a convenient  form
\begin{align}
\label{g1_2}
&\langle {\rm div}\left(S(Y_{n})-S^{\ast }(t)\right), Y_{n}-P_{n}Y\rangle\notag\\
&=\langle {\rm div} (S(Y_{n})-S(Y)), Y_{n}-P_{n}Y\rangle+\langle  {\rm div}(S(Y)-S^{\ast }(t)),Y_{n}-P_{n}Y\rangle\notag\\
&=\langle {\rm div} (S(Y_{n})-S(Y)), Y_{n}-Y\rangle
+\langle {\rm div} (S(Y_{n})-S(Y)), Y-P_{n}Y\rangle\notag\\
&+\langle  {\rm div}(S(Y)-S^{\ast }(t)),Y_{n}-P_{n}Y\rangle=g_n^1(t)+g_n^2(t)+g_n^3(t).
\end{align}
Due to  relation \eqref{29_5_S1-S2}, we have
\begin{align}
\label{g1_1}
&g_n^1(t)=-\frac{\beta}{4}\int_{\mathcal{O}}(|A_n|^2-|A|^2)^2
-\frac{\beta}{4}\int_{\mathcal{O}}(|A_n|^2+|A|^2)|A(Y_n-Y)|^2.
\end{align}
Using the inequalities \eqref{g1_2} and \eqref{g1_1}, the equation \eqref{y13} can be written as
\begin{align}
\label{31-05-10_1}
&d\left( ||Y_{n}-P_{n}Y||_{V}^{2}\right) + 4\nu
||D(Y_{n}-P_{n}Y)||_{2}^{2} dt \notag\\
&+
\frac{\beta}{2}\int_{\mathcal{O}}(|A_n|^2-|A|^2)^2dt
+\frac{\beta}{2}\int_{\mathcal{O}}(|A_n|^2+|A|^2)|A(Y_n-Y)|^2 dt\notag \\
& =2\nu\left( \Delta(P_{n}Y-Y), Y_{n}-P_{n}Y\right)  \,dt
+2\langle B(Y_{n})-B^{\ast }(t),Y_{n}-P_{n}Y\rangle \,dt\notag\\
&+2 \langle {\rm div}\left(N(Y_{n})-N^{\ast }(t)\right), Y_{n}-P_{n}Y\rangle dt
+2\left(g_n^2(t)+g_n^3(t)\right) \,dt  \notag \\
& +\sum_{i=1}^{n}|\left( \sigma(t,Y_{n})-\sigma^{\ast }(t),e_{i}\right)|
^{2}\,dt+2\left( \sigma(t,Y_{n})-\sigma^{\ast }(t),Y_{n}-P_{n}Y\right) \,d\mathcal{W}_t.
\end{align}
We also have
\begin{align}
\label{01_09}
&\langle {\rm div}\left(N(Y_{n})-N^{\ast }(t)\right), Y_{n}-P_{n}Y\rangle\notag\\
&=\langle {\rm div}\left(N(Y_{n})-N(Y)\right), Y_{n}-Y\rangle
+\langle {\rm div}\left(N(Y_{n})-N(Y)\right), Y-P_{n}Y\rangle\notag\\
&+\langle {\rm div}\left(N(Y)-N^{\ast }(t)\right), Y_{n}-P_{n}Y\rangle=h_n^1(t)+h_n^2(t)+h_n^3(t).
\end{align}
Applying Lemma \ref{Est_N1-N2<} with $3\epsilon =\frac{\beta}{8}$, we have 
\begin{align}
\label{29_3_N1-N2}
h_n^1(t)&\leq 
\frac{\beta}{8}\int_{\mathcal{O}}|A(Y_n-Y)|^2\left(|A|^2+|A_n|^2\right)
+C_1\|Y_n-Y|^2_V\notag\\
&+\frac{C}{1-\lambda}\epsilon^{\frac{\lambda-1}{\lambda+3}}\|Y_n-P_nY\|^{\frac{4(\lambda+1)}{\lambda+3}}_{H^1}\|Y\|_{H^2}^{\frac{4}{\lambda+3}}
+\frac{C}{1-\lambda}\epsilon^{\frac{\lambda-1}{\lambda+3}}\|P_nY-Y\|^{\frac{4(\lambda+1)}{\lambda+3}}_{H^1}\|Y\|_{H^2}^{\frac{4}{\lambda+3}}
\end{align}
for any $  \lambda\in]0,1[.$ Let us set  
$$
h_n^4(t)=\frac{C}{1-\lambda}\epsilon^{\frac{\lambda-1}{\lambda+3}}\|P_nY-Y\|^{\frac{4(\lambda+1)}{\lambda+3}}_{H^1}\|Y\|_{H^2}^{\frac{4}{\lambda+3}}.
$$

Proceeding analogously with the convective term, we deduce
\begin{align}
&\langle B(Y_{n}) -B^{\ast }(t),Y_{n}-P_{n}Y\rangle\notag\\
& =\langle B(Y_{n})-B(Y),Y_{n}-Y\rangle+\langle B(Y_{n})-B(Y),Y-P_{n}Y\rangle\notag\\
&+\langle B(Y)-B^{\ast }(t),Y_{n}-P_{n}Y\rangle= b_n^1(t)+ b_n^2(t)+ b_n^3(t).  \label{oi_3}
\end{align}
In addition
\begin{align}
|b_n^1(t)|&\leq C_2\|Y\|_W\|Y_n-Y\|^2_V. \label{oi_2} 
\end{align}

Denoting by $\widetilde{\sigma}_{n}$, $\widetilde{\sigma}$ and $%
\widetilde{\sigma}^{\ast }$ the solutions of the Stokes system (\ref{GS_NS}) for $%
f=\sigma(t,Y_{n})$, $f=\sigma(t,Y)$ and $f=\sigma^{\ast }(t)$, respectively, we have
\begin{equation*}
\left( \sigma(t,Y_{n})-\sigma^{\ast }(t),e_{i}\right) =(\widetilde{\sigma}_{n}-\widetilde{\sigma}%
^{\ast },e_{i})_{V},\qquad i=1,2,\dots ,n.
\end{equation*}%
Then
\begin{equation*}
\sum_{i=1}^{n}|\left( \sigma(t,Y_{n})-\sigma^{\ast }(t),e_{i}\right) |^{2}=\Vert P_{n}%
\widetilde{\sigma}_{n}-P_{n}\widetilde{\sigma}^{\ast }\Vert _{V}^{2}.
\end{equation*}%
The standard relation $x^{2}=(x-y)^{2}-y^{2}+2xy\;$ allows to write
\begin{align*}
\Vert P_{n}\widetilde{\sigma}_{n}-P_{n}\widetilde{\sigma}^{\ast }\Vert _{V}^{2}&
=\Vert P_{n}\widetilde{\sigma}_{n}-P_{n}\widetilde{\sigma}\Vert _{V}^{2}-\Vert P_{n}%
\widetilde{\sigma}-P_{n}\widetilde{\sigma}^{\ast }\Vert _{V}^{2} \\
& +2(P_{n}\widetilde{\sigma}_{n}-P_{n}\widetilde{\sigma}^{\ast },P_{n}\widetilde{\sigma}%
-P_{n}\widetilde{\sigma}^{\ast })_{V}.
\end{align*}%
From the properties of the solutions of the Stokes system (\ref{GS_NS}) and \eqref{LG_2}, we have
\begin{equation*}
\Vert P_{n}\widetilde{\sigma}_{n}-P_{n}\widetilde{\sigma}\Vert _{V}^{2}\leq \Vert
\widetilde{\sigma}_{n}-\widetilde{\sigma}\Vert _{V}^{2}\leq \Vert
\sigma(t,Y_{n})-\sigma(t,Y)\Vert _{2}^{2}\leq K\left\Vert Y_{n}-Y\right\Vert
_{V}^{2},
\end{equation*}
then
\begin{align}
\Vert P_{n}\widetilde{\sigma}_{n}-P_{n}\widetilde{\sigma}^{\ast }\Vert _{V}^{2}& \leq
K\left\Vert Y_{n}-Y\right\Vert _{V}^{2}-\Vert P_{n}\widetilde{\sigma}-P_{n}%
\widetilde{\sigma}^{\ast }\Vert _{V}^{2}  \notag \\
& +2(P_{n}\widetilde{\sigma}_{n}-P_{n}\widetilde{\sigma}^{\ast },P_{n}\widetilde{\sigma}%
-P_{n}\widetilde{\sigma}^{\ast })_{V}  \notag \\
& \leq 2K\left\Vert Y_{n}-P_{n}Y\right\Vert _{V}^{2}+C\left\Vert
P_{n}Y-Y\right\Vert _{V}^{2}-\Vert P_{n}\widetilde{\sigma}-P_{n}\widetilde{\sigma}%
^{\ast }\Vert _{V}^{2}  \notag \\
& +2(P_{n}\widetilde{\sigma}_{n}-P_{n}
\widetilde{\sigma}^{\ast },P_{n}\widetilde{\sigma}%
-P_{n}\widetilde{\sigma}^{\ast })_{V}.  \label{GG1}
\end{align}
Let us set  $D_{3}:=2(K+2C_1)$ and $D_4:=2C_2.$ 
The positive constants $K$, $C_{1}$ and  $C_{2}$ and  in \eqref{LG_2}, (\ref{29_3_N1-N2})
and  (\ref{oi_2}) 
are independent of $n$.

We introduce the auxiliary function
\begin{equation*}
\xi_2 (t)=e^{-D_{3}t-2D_{4}\int_{0}^{t}\left\Vert Y\right\Vert _{W
}ds}.
\end{equation*}
Now, applying the It\^{o} formula and  using  the  equality (\ref{31-05-10_1}), we get
\begin{align*}
& d\left( \xi_2 (t)||Y_{n}-P_{n}Y||_{V}^{2}\right) 
+4\nu\xi_2 (t)||D(Y_{n}-P_{n}Y)||_{2}^{2}\,dt \notag\\
&+ \frac{\beta}{2}\xi_2 (t)\||A_n|^2-|A|^2\|^2_2\,dt
+ \frac{\beta}{2} \xi_2 (t) \|\sqrt{|A_n|^2+|A|^2}\,|A(Y_n-Y)|\|_2^2 \,dt\notag \\
&=2\nu \xi_2 (t)(\Delta (P_{n}Y-Y),Y_{n}-P_{n}Y)\,dt  \notag \\
& +2\xi_2 (t)\langle B(Y_{n})-B^{\ast }(t),Y_{n}-P_{n}Y\rangle \,dt\notag\\
&+2\xi_2 (t)\langle {\rm div}N(Y_{n})-{\rm div}N^{\ast }(t), Y_{n}-P_{n}Y\rangle\, dt
+2\xi_2 (t)\left(g_n^2(t)+g_n^3(t)\right) \,dt   \notag \\
&+\xi_2
(t)\sum_{i=1}^{n}|\left( \sigma(t,Y_{n})-\sigma^{\ast }(t),e_{i}\right)| ^{2}\,dt  \notag\\
&+ 2\xi_2 (t)\left( \sigma(t,Y_{n})-\sigma^{\ast }(t),Y_{n}-P_{n}Y\right) \,d\mathcal{W}_t
\notag \\
& -D_{3}\xi_2 (t)||P_{n}Y-Y_{n}||_{V}^{2}\,dt-2D_{4}\xi_2 (t)\left\Vert
Y\right\Vert _{W}||Y_{n}-P_{n}Y||_{V}^{2}\,dt. 
\end{align*}
Incorporate in this equation the relations  \eqref{01_09}, \eqref{29_3_N1-N2}, \eqref{oi_3}, \eqref{oi_2} and
\eqref {GG1}, we deduce
\begin{align*}
& d\left( \xi_2 (t)||Y_{n}-P_{n}Y||_{V}^{2}\right) 
+4\nu\xi_2 (t)||D(Y_{n}-P_{n}Y)||_{2}^{2}\,dt + \frac{\beta}{2}\xi_2 (t)\||A_n|^2-|A|^2\|^2_2\,dt
\notag\\
&+ \frac{\beta}{4} \xi_2 (t) \|\sqrt{|A_n|^2+|A|^2}\,|A(Y_n-Y)|\|_2^2 \,dt +\xi_2 (t)\|P_n\widetilde{\sigma}-
P_n\widetilde{\sigma}^*\|^2_Vdt\notag\\
&\leq 2\nu \xi_2 (t)(\Delta (P_{n}Y-Y),Y_{n}-P_{n}Y)\,dt +\xi_2 (t)\frac{2C}{1-\lambda}\epsilon^{\frac{\lambda-1}{\lambda+3}}\|P_nY-Y_n\|^{\frac{4(\lambda+1)}{\lambda+3}}_{H^1}\|Y\|_{H^2}^{\frac{4}{\lambda+3}}dt\notag\\
&+2\xi_2 (t)\left[b_n^2(t)+b_n^3(t)+h_n^2(t)+h_n^3(t)+h_n^4(t)
+g_n^2(t)+g_n^3(t)\right] \,dt   \notag \\
&+\xi_2 (t)\left[ C(1+\|Y\|_W)\left\Vert
P_{n}Y-Y\right\Vert _{V}^{2}+2(P_{n}\widetilde{\sigma}_{n}-P_{n}\widetilde{\sigma}%
^{\ast },P_{n}\widetilde{\sigma}-P_{n}\widetilde{\sigma}^{\ast })_{V}\right] \,dt
 \notag\\
&+ 2\xi_2 (t)\left( \sigma(t,Y_{n})-\sigma^{\ast }(t),Y_{n}-P_{n}Y\right) \,d\mathcal{W}_t.
\end{align*}
Integrating over the time interval $(0,t\wedge \tau_M)$, $t\in [0,T]$,  and taking the expectation, we
derive
\begin{align}
\label{01_55}
& \mathbb{E} \left( \xi_2 (t\wedge \tau_M)||P_{n}Y(t\wedge \tau_M)-Y_{n}(t\wedge \tau_M)||_{V}^{2}\right) 
+4\nu\mathbb{E}\int_0^{t\wedge \tau_M}\xi_2 (s) 
||D(P_{n}Y-Y_{n})||_{2}^{2} \,ds \notag\\
&+ \frac{\beta}{2}\mathbb{E}\int_0^{t\wedge \tau_M}\xi_2 (s)\||A_n|^2-|A|^2\|^2_2\,ds
+ \frac{\beta}{4}\mathbb{E}\int_0^{t\wedge \tau_M} \xi_2 (s)
 \|\sqrt{|A_n|^2+|A|^2}\,|A(Y_n-Y)|\|^2_2 ds\notag \\
&+\mathbb{E}\int_0^{t\wedge \tau_M} \xi_2 (s)\|P_n\tilde{\sigma}-P_n\tilde{\sigma^*}\|^2_Vds\\
& \leq 2\nu \mathbb{E}\int_{0}^{t\wedge \tau_M}\xi_2 (s)(\Delta
(Y-P_{n}Y),P_{n}Y-Y_{n})\,ds \notag\\
&+2\mathbb{E}\int_0^{t\wedge \tau_M}\xi_2 (s)\left[b_n^2(s)+b_n^3(s)+h_n^2(s)+h_n^3(s)+h_n^4(s)
+g_n^2(s)+g_n^3(s)\right] \,ds \notag\\
& +\mathbb{E}\int_{0}^{t\wedge \tau_M}\xi_2 (s)\left[ C(1+M)\left\Vert
P_{n}Y-Y\right\Vert _{V}^{2}+2(P_{n}\widetilde{\sigma}_{n}-P_{n}\widetilde{\sigma}%
^{\ast },P_{n}\widetilde{\sigma}-P_{n}\widetilde{\sigma}^{\ast })_{V}\right] \,ds 
  \notag \\
&+\mathbb{E}\int_0^{t\wedge \tau_M}\xi_2 (s)\frac{2C}{1-\lambda}\epsilon^{\frac{\lambda-1}{\lambda+3}}\|P_nY-Y_n\|^{\frac{4(\lambda+1)}{\lambda+3}}_{H^1}\|Y\|_{H^2}^{\frac{4}{\lambda+3}}dt
\notag\\
& =J_n^{1}+J_n^{2}+J_n^{3}+J^4_n.
\end{align}
Here, we assume that
\begin{align}
\label{01_60}
r_n(t)=J_n^{1}+J_n^{2}+J_n^{3}\to 0.
\end{align}
This result will be proved in a lemma at the end of this  proposition.

Let us define
\begin{align}
\label{01_56}
 a_n(t)=&\mathbb{E} \left( \xi_2 (t\wedge \tau_M)||P_{n}Y(t\wedge \tau_M)-Y_{n}(t\wedge \tau_M)||_{V}^{2}\right) 
+4\nu\mathbb{E}\int_0^{t}1_{[0, \tau_M]}(s)\xi_2 (s) 
||D(P_{n}Y-Y_{n})||_{2}^{2} \,ds \notag\\
&+ \frac{\beta}{2}\mathbb{E}\int_0^{t}1_{[0, \tau_M]}(s)\xi_2 (s)\||A_n|^2-|A|^2\|^2_2\,ds\notag\\
&
+ \frac{\beta}{4}\mathbb{E}\int_0^{t}1_{[0, \tau_M]}(s) \xi_2 (s)
 \|\sqrt{|A_n|^2+|A|^2}\,|A(Y_n-Y)|\|^2_2 ds\notag \\
&+\mathbb{E}\int_0^{t}1_{[0, \tau_M]}(s) \xi_2 (s)\|P_n\widetilde{\sigma}-P_n\widetilde{\sigma}^*\|^2_Vds.
\end{align}
Taking into account \eqref{01_55},  \eqref{01_60} and  the concavity of the function
 $x\to x^{\frac{2(\lambda+1)}{\lambda+3}}$, $\lambda\in[0,1]$, we derive
\begin{align}
\label{01_57}
a_n(t)&\leq r_n(t)+\frac{2C}{1-\lambda}\mathbb{E}\int_0^{t}1_{[0, \tau_M]}(s)\xi_2 (s)\epsilon^{\frac{\lambda-1}{\lambda+3}}\|P_nY-Y_n\|^{\frac{4(\lambda+1)}{\lambda+3}}_{H^1}\|Y\|_{H^2}^{\frac{4}{\lambda+3}}ds\notag\\
&\leq r_n(t)+ M^{\frac{4}{\lambda+3}} \frac{2C}{1-\lambda}
\epsilon^{\frac{\lambda-1}{\lambda+3}}
\mathbb{E}\int_0^{t}
\xi_2 (s\wedge \tau_M)\|P_nY(s\wedge \tau_M)-Y_n(s\wedge \tau_M)\|^{\frac{4(\lambda+1)}{\lambda+3}}_{H^1}ds
\notag\\
&\leq r_n(t)+ M^{\frac{4}{\lambda+3}} \frac{2C}{1-\lambda}
\epsilon^{\frac{\lambda-1}{\lambda+3}}
\mathbb{E}\int_0^{t}
\left(\xi_2 (s\wedge \tau_M)\|P_nY(s\wedge \tau_M)-Y_n(s\wedge \tau_M)\|^2_V\right)^{\frac{2(\lambda+1)}{\lambda+3}}ds
\notag\\
&\leq r_n(t)+ M^{\frac{4}{\lambda+3}} \frac{2C}{1-\lambda}
\epsilon^{\frac{\lambda-1}{\lambda+3}}
\left(\int_0^{t}\mathbb{E}\,\xi_2 (s\wedge \tau_M)\|P_nY(s\wedge \tau_M)-Y_n(s\wedge \tau_M)\|^2_Vds\right)^{\frac{2(\lambda+1)}{\lambda+3}}
\notag\\&\leq  r_n(t)+M^{\frac{4}{\lambda+3}} \frac{2C}{1-\lambda}
\epsilon^{\frac{\lambda-1}{\lambda+3}}
\left(\int_0^{t}a_n(s)ds\right)^{\frac{2(\lambda+1)}{\lambda+3}},
\end{align}
which yields 
\begin{align}
\label{10:26_07:08}
\limsup_{n\to\infty} a_n(t)
\leq \limsup_{n\to\infty} r_n(t)+M^{\frac{4}{\lambda+3}} \frac{2C}{1-\lambda}
\epsilon^{\frac{\lambda-1}{\lambda+3}}
\left(\int_0^{t}\limsup_{n\to\infty}a_n(s)ds\right)^{\frac{2(\lambda+1)}{\lambda+3}}.
\end{align}
Denoting
$$
f(t):=\int_0^{t}\limsup_{n\to\infty}\,a_n(s)ds,
$$
and knowing that $\lim_{n\to\infty} r_n(t)=0$, \eqref{10:26_07:08} can be written as 
\begin{align}
\label{01_10_42}
f'(t)
&\leq  M^{\frac{4}{\lambda+3}} \frac{2C}{1-\lambda}
\epsilon^{\frac{\lambda-1}{\lambda+3}}
\left(f(t)\right)^{\frac{4(\lambda+1)}{\lambda+3}}.
\end{align}
Here, we can proceed as in \cite{BI04} in order to verify that $f\equiv 0.$
Since $f(0)=0$ and 
$$
\left((f(t))^{\frac{1-\lambda}{\lambda+3}}\right)^\prime\leq \frac{2C}{\lambda+3}
M^{\frac{4}{\lambda+3}}\epsilon^{\frac{\lambda+1}{\lambda+3}},
$$
we have
$$
f(t)\leq \left(\frac{2C}{\lambda+3}
M^{\frac{4}{\lambda+3}}\epsilon^{\frac{\lambda+1}{\lambda+3}}\,t
\right)^{^{\frac{\lambda+3}{1-\lambda}}}.
$$
Considering  $T_0=\frac{3}{4CM^{\frac{4}{3}}\epsilon}$, we have  $\frac{2C}{\lambda+3}
M^{\frac{4}{\lambda+3}}\epsilon^{\frac{\lambda+1}{\lambda+3}}\,t\leq \frac{1}{2}$. Taking 
$\lambda\to 1$, we get $f(t)=0,\, \forall t\in [0,T_0].$ By an extension argument, we obtain   $f(t)=0,\, \forall t\in [0,T].$

\begin{align*}
& \mathbb{E} \left( \xi_2 (t\wedge \tau_M)||P_{n}Y(t\wedge \tau_M)-Y_{n}(t\wedge \tau_M)||_{V}^{2}\right) 
+ 4\nu\mathbb{E}\int_0^{t\wedge \tau_M}\xi_2 (s)
||D(P_{n}Y-Y_{n})||_{2}^{2}\, ds \notag\\
&+ 
\frac{\beta}{2}\mathbb{E}\int_0^{t\wedge \tau_M}\xi_2 (s)\int_{\mathcal{O}}(|A_n|^2-|A|^2)^2
ds
+ \frac{\beta}{4}\mathbb{E}\int_0^{t\wedge \tau_M} \xi_2 (s) \int_{\mathcal{O}}(|A_n|^2+|A|^2)|A(Y_n-Y)|^2 ds\notag \\
&+\mathbb{E}\int_0^{t\wedge \tau_M}
 \xi_2 (s)\|P_n\widetilde{\sigma}-P_n\widetilde{\sigma}\|^2_Vds \to 0, \qquad\text{as}\quad n\to\infty.
\end{align*}

\bigskip\ $\hfill \hfill \blacksquare $

\begin{lemma}
Let $J_n^1(t)$,  $J_n^2(t)$, $J_n^3(t)$ be the terms introduced  in \eqref{01_55}. Then 
for all $t\in[0,T]$, $J_n^i(t)\to 0$, for $i=1,2,3.$
\end{lemma}
\textbf{Proof.} Using (\ref{as})-(\ref{c1}) and the properties of the projection $P_{n},$ we
have
\begin{align*}
|J_n^{1}(t)|& =\bigl|2\nu \mathbb{E}\int_{0}^{t}\xi_2 (s)(1_{[0,\tau
_{M}]}(s)\Delta (Y-P_{n}Y),P_{n}Y-Y_{n})\,ds\bigr| \\
& \leq C\Vert P_{n}Y-Y\Vert _{L^{2}(\Omega \times (0,t),H^{2})}\Vert
P_{n}Y-Y_{n}\Vert _{L^{2}(\Omega \times (0,t),W)} \\
& \leq C\Vert P_{n}Y-Y\Vert _{L^{2}(\Omega \times (0,T),W)}\left( \Vert
Y\Vert _{L^{2}(\Omega \times (0,T),W)}+\Vert Y_{n}\Vert _{L^{2}(\Omega
\times (0,T);H^{2})}\right) \\
& \leq C\Vert P_{n}Y-Y\Vert _{L^{2}(\Omega \times (0,T),H^{2})},
\end{align*}
which goes to zero, as $n\rightarrow \infty $, by (\ref{c02Y}). 
$$
J_n^2(t)=
2\mathbb{E}\int_0^{t\wedge \tau_M}\xi_2 (s)\left[b_n^2(s)+b_n^3(s)+h_n^2(s)+h_n^3(s)+h_n^4(s)
+g_n^2(s)+g_n^3(s)\right] \,ds
$$
From \eqref{14:31_07:08}, (\ref{ineq2}) and \eqref{c02Y}, we deduce
\begin{align*}
\left|2\mathbb{E}\int_0^{t\wedge \tau_M}\xi_2 (s)b_n^2(s)\right|\, ds
&=\left|2\mathbb{E}\int_0^{t\wedge \tau_M}\xi_2 (s)\langle B(Y_{n})-B(Y),P_{n}Y-Y\rangle\right|\notag\\
&\leq C\mathbb{E}\sup_{t\in[0,T]}\|Y_n\|^4_{V}\;\mathbb{E}\|P_{n}Y-Y\|^2_{L^2(0,T;V)}
\to 0, \quad \text{as}\quad n\to \infty.
\end{align*}
Convergences (\ref{c1}) and (\ref{c02Y}) give that
\begin{equation*}
P_{n}Y-Y_{n}\rightarrow 0\qquad \text{ weakly in }L^{2}(\Omega \times (0,T),%
W),
\end{equation*}%
then for any operator $R\in L^{2}(\Omega \times (0,T),W^{\ast })$
we have
\begin{equation*}
\mathbb{E}\int_{0}^{T}\langle R,P_{n}Y-Y_{n}\rangle \,ds\rightarrow 0,\qquad
\text{as }n\rightarrow \infty .
\end{equation*}
The function $1_{[0,\tau _{M}]}(s)\xi_2 (s)$ is bounded, then
\begin{eqnarray*}
&&\Vert 1_{[0,\tau _{M}]}(s)\xi_2 (s)\left( B(Y)-B^{\ast }\right) \Vert
_{L^{2}(\Omega \times (0,T),W^{\ast })}^{2} \\
&\leq &C\left( \Vert B(Y)\Vert _{L^{2}(\Omega \times (0,T),W
^{\ast })}^{2}+\Vert B^{\ast }\Vert _{L^{2}(\Omega \times (0,T), W
^{\ast })}^{2}\right) \leq C,
\end{eqnarray*}%
by (\ref{as}), (\ref{l}) and (\ref{c01}). Therefore, $\text{as }n\rightarrow \infty$, we have
\begin{align*}
2\mathbb{E}\int_0^{t\wedge \tau_M}\xi_2 (s)b_n^3(s) & =2\mathbb{E}\int_{0}^{t}\langle 1_{[0,\tau _{M}]}(s)\xi_2 (s)\left(
B(Y)-B^{\ast }(s)\right) ,P_{n}Y-Y_{n}\rangle \,ds\rightarrow 0.
\end{align*}
Using the same reasoning, we show that
\begin{align*}
2\mathbb{E}\int_0^{t\wedge \tau_M}\xi_2 (s)h_n^2(s)\to 0,\quad 
2\mathbb{E}\int_0^{t\wedge \tau_M}\xi_2 (s)h_n^3(s)\to 0.
\end{align*}
By the definition of the stopping time $\tau_M$, we have $1_{[0,\tau_M]}(s)\xi_2 (s)\|Y\|_{W}^{\frac{4}{\lambda+3}}\leq M^{\frac{4}{\lambda+3}}$, so
\begin{align*}
\left|2\mathbb{E}\int_0^{t\wedge \tau_M}\xi_2 (s)h_n^4(s)\right| 
&\leq \frac{C}{1-\lambda}\epsilon^{\frac{\lambda-1}{\lambda+3}} \left|2\mathbb{E}\int_0^{T}1_{[0,\tau_M]}(s)\xi_2 (s)\|P_nY-Y\|^{\frac{4(\lambda+1)}{\lambda+3}}_{V}\|Y\|_{W}^{\frac{4}{\lambda+3}}\right|\notag\\
&\leq \frac{C}{1-\lambda}\epsilon^{\frac{\lambda-1}{\lambda+3}}
 M^{\frac{4}{\lambda+3}}
 \left|2\mathbb{E}
\int_0^{T}\|P_nY-Y\|^{\frac{4(\lambda+1)}{\lambda+3}}_{V}
\right|\notag\\
&\leq C(M,\lambda)
\,
\|P_nY-Y\|^{\frac{2(\lambda+1)}{\lambda+3}}_{L^2(\Omega\times(0,T), Y)}
\to 0
\end{align*}
Similarly we verify that the remaining terms in $J_n^2(t)$ converges to $0$,
as well as  $J_n^3(t)$ converges to $0$, as $n\to \infty$.

\bigskip\ $\hfill \hfill \blacksquare $

From \eqref{15:47_07:08},  the following strong convergences hold 
\begin{align}
\label{Lim_1}
\lim_{n\rightarrow \infty }\mathbb{E}\left( \xi_2 (\tau _{M})||P_{n}Y(\tau
_{M})-Y_{n}(\tau _{M})||_{V}^{2}\right) & =0, \\
\label{Lim_2}
\lim_{n\rightarrow \infty }\mathbb{E}\int_{0}^{\tau _{M}}\xi_2 (s)
||D(P_{n}Y-Y_{n})||_{2}^{2} ds& =0,\\
\label{Lim_3}
\lim_{n\rightarrow \infty }\mathbb{E}\int_0^{ \tau_M}\xi_2 (s)\||A_n|^2-|A|^2\|^2_2\,ds&=0,\\
\label{Lim_4}
\lim_{n\rightarrow \infty }\mathbb{E}\int_0^{ \tau_M} \xi_2 (s)
 \|\sqrt{|A_n|^2+|A|^2}\,|A(Y_n-Y)|\|^2_2 ds&=0,\\
 \label{Lim_5}
 \lim_{n\rightarrow \infty }\mathbb{E}\int_{0}^{\tau _{M}}\xi_2 (s)\Vert P_{n}%
\widetilde{\sigma}-P_{n}\widetilde{\sigma}^{\ast }\Vert _{V}^{2}ds& =0,
\end{align}
for each $M\in \mathbb{N}$. Since there exists a strictly positive constant $%
\mu,$ such that $\mu\leq 1_{[0,\tau _{M}]}(s)\xi_2 (s)\leq
1,$ it follows from the Korn inequality \eqref{14:37_3Ag} and (\ref{c02Y}) that
\begin{equation}
\label{13-16-34}
\lim_{n\rightarrow \infty }\mathbb{E}\int_{0}^{\tau _{M}}\xi_2 (s)
||D(P_{n}Y-Y_{n})||_{2}^{2} ds=0
\quad \text{implies}\quad
\lim_{n\rightarrow \infty }\mathbb{E}\int_{0}^{\tau
_{M}}||Y-Y_{n}||_{V}^{2}ds=0.
\end{equation}
In addition, we have
\begin{align}
\label{Lim_3}
\lim_{n\rightarrow \infty }\mathbb{E}\int_0^{ \tau_M}\||A_n|^2-|A|^2\|^2_2\,ds&=0,\\
\label{Lim_4}
\lim_{n\rightarrow \infty }\mathbb{E}\int_0^{ \tau_M} 
 \|\sqrt{|A_n|^2+|A|^2}\,|A(Y_n-Y)|\|^2_2 ds&=0.
\end{align}
Considering (\ref{c02}), we also derive
\begin{equation}
\mathbb{E}\int_{0}^{\tau _{M}}\Vert \widetilde{\sigma}-\widetilde{\sigma}^{\ast }\Vert
_{V}^{2}ds=0.  \label{RC2}
\end{equation}

\bigskip

\textit{Step 4. Identification of }
$B^{\ast }(t)$\textit{ with} $B(Y)$, ${\rm div}\,N^{\ast }(t)$\textit{ with} ${\rm div}\, N(Y)$,
${\rm div}\,S^{\ast }(t)$\textit{ with} ${\rm div}\, S(Y)$
\textit{\ and }
$\sigma^{\ast }(t)$ \textit{ with} $\sigma(t,Y)$ \textit{on} $[0,\tau_M]$ \textit{for each} $M.$

\bigskip

Now, we are able to show that the limit function $Y$ satisfies equation (\ref%
{var_form_state}). Integrating equation (\ref{y12}) on the time interval $%
(0,\tau _{M}\wedge t)$, we derive
\begin{eqnarray}
\left( \upsilon \left( Y(\tau _{M}\wedge t)\right) ,\phi \right) -\left(
\upsilon \left( Y_{0}\right) ,\phi \right) &=&\int_{0}^{\tau _{M}\wedge t}
\bigl[ \left( \nu \Delta Y+U,\phi \right) +\langle B^{\ast }(s),\phi \rangle+\langle {\rm div}\,N^{\ast }(s),\phi \rangle \notag \\
&&+\langle {\rm div}\,S^{\ast }(s),\phi \rangle %
\bigr] \,ds +\int_{0}^{\tau _{M}\wedge t}\left( \sigma^{\ast }(s),\phi \right) \,d\mathcal{W}_s
\label{y122}
\end{eqnarray}%
for any $\phi \in V.$
From (\ref{RC2}) it follows that%
\begin{equation*}
1_{[0,\tau _{M}]}(t)\widetilde{\sigma}=1_{[0,\tau _{M}]}(t)\widetilde{\sigma}^{\ast
}\qquad \text{a.e. in }\Omega \times (0,T),
\end{equation*}%
which implies
\begin{equation}
1_{[0,\tau _{M}]}(t)\sigma(t,Y)=1_{[0,\tau _{M}]}(t)\sigma^{\ast }(t)\qquad \text{a.
e. in }\Omega \times (0,T)  \label{09}
\end{equation}%
by \ (\ref{GS_NS}). Since 
$B(Y_{n})-B(Y)=(Y_{n}\cdot \nabla)(Y_{n}-Y)+(Y_{n}-Y)\cdot \nabla Y,$ we verify that
\begin{equation*}
\Vert B(Y_{n})-B(Y)\Vert _{V^{\ast }}\leq C\left( \Vert
Y_{n}\Vert _{V}+\Vert Y\Vert _{V}\right) \Vert
Y_{n}-Y\Vert _{V}.
\end{equation*}%
Then for any $\varphi \in L^{\infty }(\Omega \times (0,T),V),$
using (\ref{as}), (\ref{c1})
\begin{align*}
\bigl\vert\mathbb{E}\int_{0}^{T}1_{[0,\tau _{M}]}(s)\langle B(Y_{n})&
-B(Y),\ \varphi \rangle \,ds\bigr\vert \\
& \leq C\mathbb{E}\int_{0}^{T}1_{[0,\tau _{M}]}(s)\left( \Vert Y_{n}\Vert _{%
V}+\Vert Y\Vert _{V}\right) \Vert Y_{n}-Y\Vert
_{V}\Vert \varphi \Vert _{V}\,ds \\
& \leq C\Vert \varphi \Vert _{L^{\infty }(\Omega \times (0,T),V)}%
\mathbb{E}\int_{0}^{T}\left( \Vert Y_{n}\Vert _{V}+\Vert Y\Vert
_{V}\right) \Vert Y_{n}-Y\Vert _{V}\,ds \\
& \leq C\Vert \varphi \Vert _{L^{\infty }(\Omega \times (0,T),V%
)}\left( \mathbb{E}\int_{0}^{\tau _{M}}\Vert Y_{n}-Y\Vert
_{V}^{2}\,ds\right) ^{\frac{1}{2}}\rightarrow 0,\qquad \text{as }%
n\rightarrow \infty .
\end{align*}%
Taking into account  \ (\ref{c01})$_{1}$ and that the space $L^{\infty }(\Omega \times
(0,T),V)$ is dense in $L^{2}(\Omega \times (0,T),V),$
we obtain
\begin{equation}
1_{[0,\tau _{M}]}(s)B^{\ast }(s)=1_{[0,\tau _{M}]}(s)B(Y)\qquad \text{a.
e. in }\Omega \times (0,T).  \label{RC2Y}
\end{equation}

From (\ref{13-15-22}), we have
\begin{align}
\label{13-16-14}
\left|\langle{\rm div}(N({Y_n})- N(Y)),\phi\rangle\right|
&\leq 
C\epsilon \left\|A(Y_n-Y)\sqrt{|A_n|^2+|A|^2}\right\|_2\|\phi\|_V\notag\\
&
+C\|Y_n-y\|_V\left(\|Y_n\|_W+\|Y\|_W\right)\|\phi\|_W.
\end{align}
Then for any $\phi \in L^{\infty }(\Omega \times (0,T),W),$
using (\ref{13-16-34}) and  (\ref{Lim_4}), we deduce
\begin{align*}
&\left|\mathbb{E}\int_{0}^{T}1_{[0,\tau _{M}]}(s)\langle{\rm div}(N(Y_n)- N(Y)),\phi\rangle \,ds\right| \\
& \leq C\mathbb{E}\int_{0}^{T}1_{[0,\tau _{M}]}(s)
 \left\|A(Y_n-Y)\sqrt{|A_n|^2+|A|^2}\right\|_2\|\phi\|_V
\,ds \\
& +  C\mathbb{E}\int_{0}^{T}1_{[0,\tau _{M}]}(s)
\|Y_n-y\|_V\left(\|Y_n\|_W+\|Y\|_W\right)\|\phi\|_W\notag\\
&\leq C\Vert \phi \Vert _{L^{\infty }(\Omega \times (0,T),W)}
\mathbb{E}\int_{0}^{\tau_M}
 \left\|A(Y_n-Y)\sqrt{|A_n|^2+|A|^2}\right\|_2
\,ds \\
&+C\Vert \phi \Vert _{L^{\infty }(\Omega \times (0,T),W)}
\mathbb{E}\int_{0}^{\tau_M}\left( \Vert Y_{n}\Vert _{W}+\Vert Y\Vert
_{W}\right) \Vert Y_{n}-Y\Vert _{V}\,ds \\
& \leq
 C\Vert \phi \Vert _{L^{\infty }(\Omega \times (0,T),W)}
\mathbb{E}\int_{0}^{\tau_M}
 \left\|A(Y_n-Y)\sqrt{|A_n|^2+|A|^2}\right\|^2_2
\,ds \\
&+ C\Vert \phi \Vert _{L^{\infty }(\Omega \times (0,T),W
)}\left( \mathbb{E}\int_{0}^{\tau _{M}}\Vert Y_{n}-Y\Vert
_{V}^{2}\,ds\right) ^{\frac{1}{2}}\rightarrow 0,\qquad \text{as }%
n\rightarrow \infty .
\end{align*}
Therefore
\begin{equation}
1_{[0,\tau _{M}]}(s)\,{\rm div}(N^\ast (s) )=1_{[0,\tau _{M}]}(s)\,{\rm div}(N(Y) )\qquad \text{a.e. in }\Omega \times (0,T).  \label{13-16-41}
\end{equation}
Using the same reasoning, we  show
\begin{equation}
1_{[0,\tau _{M}]}(s)\,{\rm div}(S^\ast (s) )=1_{[0,\tau _{M}]}(s)\,{\rm div}(S(Y) )\qquad \text{a.e. in }\Omega \times (0,T).  \label{13-16-42}
\end{equation}
Namely, from \eqref{29_5_S1-S2}, we have
\begin{align}
\left|\langle {\rm div} (S(Y)-S(Y_n)), \phi\rangle\right|\leq 
C\|Y\|^2_{W}\|Y-Y_n\|_{V}\|\phi\|_{W}+
C\|Y_n\|_{W}\left\||A|^2-|A_n|^2\right\|_2\|\phi\|_{W},
\end{align}
and (\ref{13-16-34}) and  (\ref{Lim_3}) gives
\begin{align*}
&\left|\mathbb{E}\int_{0}^{T}1_{[0,\tau _{M}]}(s)
\langle {\rm div} (S(Y)-S(Y_n), \phi\rangle \,ds\right| \\
& \leq C\mathbb{E}\int_{0}^{T}1_{[0,\tau _{M}]}(s)
 C\|Y\|^2_{W}\|Y-Y_n\|_{V}\|\phi\|_{W}
\,ds \\
& +  C\mathbb{E}\int_{0}^{T}1_{[0,\tau _{M}]}(s)\|Y_n\|_{W}\left\||A|^2-|A_n|^2\right\|_2\|\phi\|_{W}\notag\\
&\leq C\Vert \phi \Vert _{L^{\infty }(\Omega \times (0,T),W)}
\mathbb{E}\int_{0}^{\tau_M}\Vert Y\Vert
_{W} \Vert Y_{n}-Y\Vert _{V}\,ds \\
&+C\Vert \phi \Vert _{L^{\infty }(\Omega \times (0,T),W)}
\mathbb{E}\int_{0}^{\tau_M}\|Y_n\|_{W}\left\||A|^2-|A_n|^2\right\|_2\\
&\leq C(M)\Vert \phi \Vert _{L^{\infty }(\Omega \times (0,T),W)}
\mathbb{E}\int_{0}^{\tau_M}\Vert Y_{n}-Y\Vert _{V}\,ds \\
&+C(M)\Vert \phi \Vert _{L^{\infty }(\Omega \times (0,T),W)}
\mathbb{E}\int_{0}^{\tau_M}\left\||A|^2-|A_n|^2\right\|_2 \rightarrow 0,\qquad \text{as }%
n\rightarrow \infty .
\end{align*}

By introducing identities (\ref{09}), (\ref{RC2Y}), \eqref{13-16-41} and 
\eqref{13-16-42} in equation (\ref{y122}),
it follows that
\begin{align}
\left( \upsilon \left( Y(\tau _{M}\wedge t)\right) ,\phi \right) &-\left(
\upsilon \left( Y_{0}\right) ,\phi \right)  =\int_{0}^{\tau _{M}\wedge t}
\biggl[ \left( \nu \Delta Y+U,\phi \right) +\langle B(Y)\notag\\
&+{\rm div} (N(Y))+{\rm div} (S(Y)),\phi \rangle 
\biggr] \,ds+\int_{0}^{\tau _{M}\wedge t}\left( \sigma(s,Y),\phi \right) \,d\mathcal{W}_s.
\label{y1122}
\end{align}
Reasoning as in (\ref{fc})  we have $\tau _{M}\rightarrow T$
 \ a.e. in
$\Omega $, as  $M\to\infty$. We can pass to the limit in each term of equation (\ref{y1122})
in $L^{1}(\Omega \times (0,T))$, as $M\rightarrow \infty, $  by applying the
Lebesgue dominated convergence theorem and the Burkholder-Davis-Gundy
inequality for the last (stochastic) term, deriving an equivalent formulation of 
equation (\ref{var_form_state}) a.e. in $\Omega \times (0,T)$.

\bigskip

\textit{Step 5. Uniqueness.} In order to prove uniqueness,
 we take two solutions $Y_1$ and $Y_1$, and consider the difference $Y=Y_1-Y_2$. 
Using  similar arguments as in the previous steps, introducing the function  
\begin{equation*}
\xi_3 (t)=e^{-\frac{1}{2}D_{3}t-D_{4}\int_{0}^{t}\left\Vert Y_1\right\Vert _{W
}ds},
\end{equation*}
 we  show that 
\begin{align*}
 \mathbb{E} \left( \xi_3 (t)||Y(t)||_{V}^{2}\right)
&+ 4\nu\mathbb{E}\int_0^{t}\xi_3 (s)
||D(Y)||_{2}^{2} ds \notag\\
&+ 
\frac{\beta}{2}\mathbb{E}\int_0^{t}\xi_3 (s)\int_{\mathcal{O}}(|A(Y_1)|^2-|A(Y_2)|^2)^2ds\notag\\
&
+ \frac{\beta}{4}\mathbb{E}\int_0^{t} \xi_3 (s) \int_{\mathcal{O}}(|A(Y_1)|^2+|A(Y_2)|^2)|A(Y)|^2 ds\notag \\
&=0\quad \text{for a.e. } t\in[0,T].
\end{align*}
Therefore, for a.e. $t\in[0,T]$, we have 
$$
\mathbb{E} \left( \xi_3 (t)||Y(t)||_{V}^{2}\right)=0.
$$
Since $\xi_3$ is a positive function, we deduce that 
for a.e. $t\in[0,T]$
$$
\quad Y_1(t)=Y_2(t), \quad P-\text{a.s.}.
$$

 $\hfill \hfill \blacksquare $

\bigskip

\appendix
\section{Appendix} \label{appendix}
\setcounter{equation}{0}
In this appendix, we collect  important inequalities from \cite{BI06} related with the nonlinear 
terms that we apply throughout the article.
  \begin{lemma}
\label{09-17-29}
For any $\epsilon >0$ and  $y\in W$, we have
\begin{align}
\label{08-11-48}
\left|
(\alpha_1+\alpha_2)\int_{\mathcal{O}}
{\rm div}
(A^2)\cdot y\right|&\leq \epsilon \|A^2\|^2_2+
\frac{(\alpha_1+\alpha_2)^2}{16\epsilon}
\|A\|^2_2,
\end{align}
where $A=A(y).$
\end{lemma}
\textbf{Proof.}
The integration by parts gives
\begin{align}
(\alpha_1+\alpha_2)\int_{\mathcal{O}}
{\rm div}
(A^2)\cdot y=
(\alpha_1+\alpha_2)\int_{\Gamma}
\left(\mathrm{n}\cdot A^2\right)\cdot  y-
(\alpha_1+\alpha_2)\int_{\mathcal{O}}
A^2\cdot \nabla y.
\end{align}
Due to the boundary conditions $y=(y\cdot \tau)\tau$ and $(\mathrm{n}\cdot A)\cdot \tau=0$ on $\Gamma$,   we obtain
\begin{align*}
(\mathrm{n}\cdot  A^2)\cdot y&=(y\cdot \tau)(\mathrm{n}\cdot  A^2)\cdot\tau=(y\cdot \tau)(\mathrm{n}\cdot  A^2)\cdot \tau=(y\cdot \tau)((\mathrm{n}\cdot  A)\cdot A)\cdot \tau\\
&=(y\cdot \tau)\left[((\mathrm{n}\cdot  A )\cdot\mathrm{n} )((\mathrm{n}\cdot A)\cdot \tau)
+( (\mathrm{n}\cdot  A)\cdot \mathrm{\tau })((\mathrm{\tau }\cdot A)\cdot \tau)\right]=0.
\end{align*}
Using the symmetry of $A$, we derive
\begin{align}
(\alpha_1+\alpha_2)\int_{\mathcal{O}}
{\rm div}
(A^2)\cdot y=-
\frac{1}{2}(\alpha_1+\alpha_2)\int_{\mathcal{O}}
A^2\cdot A.
\end{align}
Therefore, the H\"older and the  Young inequalities give 
\eqref{08-11-48}.

\bigskip\ $\hfill \hfill \blacksquare $

Considering a small change in  estimate $(35)$ of \cite{BI06}, we collect the following  estimates.
\begin{lemma}[see \cite{BI06}, relations (33)-(36)]
\label{Est_K} For each $y\in W$ and any $ \epsilon, \,\delta>0$, the following estimates are valid
\begin{align}
\label{09-13-02}
\left({\rm div}\left(|A|^2A\right),\mathbb{P}\upsilon (y)\right)&\leq -\frac{1}{2}\|A\|_4^4-\frac{\alpha_1}{2}\||A||\nabla A|\|^2_2-\frac{\alpha_1}{4}\|\nabla(|A|^2)\|^2_2\notag\\
&\qquad +3\epsilon \||A||\nabla^2y|\|^2_2+5\epsilon\|A\|^4_{12}+3\epsilon\|y\|^4_{H^1}+C(\epsilon)
\|y\|_{H^1}^2\|y\|^2_{H^2},
\end{align}
\begin{align}
\label{09-13-03}
(\alpha_1+\alpha_2)\bigl({\rm div}\left(A^2\right),
\mathbb{P}\upsilon (y)\bigr)&\leq 
 \epsilon \||A||\nabla^2y|\|^2_2+C(\epsilon)\|y\|^2_W,
 \end{align}
 \begin{align}
\label{09-13-04}
 -\bigl((y\cdot \nabla) \upsilon+\sum_j\upsilon^j  \nabla y^j,\mathbb{P}\upsilon (y)\bigr)&\leq 4\epsilon \||A||\nabla^2y|\|^2_2
 +C(\epsilon, \delta)\|y\|^2_{W}
 +C(\epsilon)\|y\|_\infty\|y\|^2_{W}
 +\delta \|y\|^4_{W^{1,4}},
\end{align}
where $A=A(y)$ and $\upsilon=\upsilon(y).$ 

\end{lemma}

Let us introduce the operators
\begin{eqnarray}
&&S(y):=\beta\left(|A(y)|^2A(y)\right), \quad \label{operator1}\\
&&N(y):=\alpha_1\left(y\cdot \nabla A(y)+(\nabla y)^\top\, A(y)+A(y)\, \nabla y\right)
-\alpha_2(A(y))^2.\label{operator2}
\end{eqnarray}

\begin{lemma} For any $y, \,\hat{y}\in W$, we have
\begin{align}
\label{29_5_S1-S2}
\langle {\rm div} (S(\hat{y})-S(y)), \hat{y}-y\rangle= 
-\frac{\beta}{4}\int_{\mathcal{O}}(|\hat{A}|^2-|A|^2)^2
-\frac{\beta}{4}\int_{\mathcal{O}}(|\hat{A}|^2+|A|^2)|A(\hat{y}-y)|^2,
\end{align}
where $A=A(y)$ and $\hat{A}=A(\hat{y}).$
\end{lemma}
\textbf{Proof.}  Integrating  by parts, we write
\begin{align}
\langle {\rm div} (S(\hat{y})-S(y)), \hat{y}-y\rangle
&=\beta\int_{\Gamma}
\left(\mathrm{n}\cdot (|\hat{A}|^2\hat{A}-|A|^2A)\right)\cdot (\hat{y}-y)\notag\\
&-\frac{\beta}{2}\int_{\mathcal{O}}(|\hat{A}|^2\hat{A}-|A|^2A)\cdot A(\hat{y}-y)
=I_{11}+I_{12}.
\end{align}
Using the boundary conditions, we deduce that
\begin{align}
I_{11}&=\beta\int_{\Gamma}((\hat{y}-y)\cdot \tau)
\left[|\hat{A}|^2(\mathrm{n}\cdot \hat{A})\cdot \tau
-|A|^2(\mathrm{n}\cdot A)\cdot \tau
\right]=0.
\end{align}
Standard algebraic computations yield
\begin{align}
I_{12}=
&-\frac{\beta}{2}\int_{\mathcal{O}}(|\hat{A}|^2\hat{A}-|A|^2A)\cdot A(\hat{y}-y)
\notag\\
&=-\frac{\beta}{2}\int_{\mathcal{O}}(|\hat{A}|^2\hat{A}-|A|^2A)\cdot A(\hat{y}-y)
\notag\\
&=-\frac{\beta}{4}\int_{\mathcal{O}}(|\hat{A}|^2-|A|^2)^2
-\frac{\beta}{4}\int_{\mathcal{O}}(|\hat{A}|^2+|A|^2)|A(\hat{y}-y)|^2.
\end{align}

\bigskip\ $\hfill \hfill \blacksquare $

\begin{lemma}
\label{Est_N1-N2<} 
For any $y,\; \hat{y} \in W$, the following estimate holds
\begin{align}
\label{25_2_N1-N2}
\left({\rm div}\left(N(\hat{y})- N(y)\right),\hat{y}-y\right)
&=\left(N(\hat{y})- N(y),\nabla(\hat{y}-y)\right)\notag\\
&\leq 
3\epsilon \int_{\mathcal{O}}|A(\hat{y}-y)|^2\left(|A|^2+|\hat{A}|^2\right)
+\frac{C}{\epsilon}\int_{\mathcal{O}}|\nabla(\hat{y}-y)|^2\notag\\
&+\frac{C}{1-\lambda}\epsilon^{\frac{\lambda-1}{\lambda+3}}\|\hat{y}-y\|^{\frac{4(\lambda+1)}{\lambda+3}}_{H^1}\|y\|_{H^2}^{\frac{4}{\lambda+3}}\quad \text{for any} \quad\epsilon>0,
\;  \lambda\in]0,1[,
\end{align}
where $A=A(y)$ and $\hat{A}=A(\hat{y})$.
\end{lemma}

\textbf{Proof.}
The divergence theorem gives
\begin{align}
\left({\rm div}\left(N(\hat{y})- N(y)\right),\hat{y}-y\right)&
=\left(N(\hat{y})- N(y),\nabla(\hat{y}-y)\right)-\int_\Gamma
\left[\left(N(\hat{y})- N(y)\right)\mathrm{n}\right]\cdot(\hat{y}-y).
\end{align}
The relation \eqref{25_2_N1-N2} is proved in  \cite{BI04} for the case $\mathcal{O}=\mathbb{R}^2$ (domain without boundary), in \cite{BI06} it is verified that 
the boundary term vanishes. 

\bigskip\ $\hfill \hfill \blacksquare $


\textbf{Acknowledgment}
The work of F. Cipriano was partially supported by the Funda\c{c}\~{a}o para
a Ci\^{e}ncia e a Tecnologia (Portuguese Foundation for Science and
Technology) through the project UID/MAT/00297/2019 (Centro de Matem\'{a}tica
e Aplica\c{c}\~{o}es).

\textbf{\bigskip }

\end{document}